\documentclass[12pt]{article}
{\begin{Sbox}\begin{minipage}}%
		{\end{minipage}\end{Sbox}\fbox{\TheSbox}}%
\usepackage[psamsfonts]{amssymb}
\usepackage{amsmath, amsthm, anysize, enumerate, color, url}
\allowdisplaybreaks[4]
\usepackage{graphicx}
\usepackage{epstopdf}
\usepackage{epsfig}
\usepackage{subfigure}
\usepackage{multirow}
\usepackage[ruled,linesnumbered]{algorithm2e}
\usepackage{tikz}
\usepackage{array}

\usepackage[sectionbib]{natbib}

\usepackage{threeparttable}

\usepackage[colorlinks, breaklinks,
linkcolor=red,
citecolor=blue
]{hyperref}
\usepackage{breakurl}
\usepackage{bm}

\newtheorem{theorem}{Theorem} 
\newtheorem{cor}{Corollary}

\newtheorem{lemma}{Lemma} 

\newcommand{\balign}{\begin{align}}
	\newcommand{\ealign}{\end{align}}
\newcommand{\beq}{\begin{equation}}
	\newcommand{\eeq}{\end{equation}}

\newtheorem{definition}{Definition}

\theoremstyle{definition}

\newcommand{\E}{\mathbb{E}}

\newcommand{\diag}{\mbox{diag}}

\newcommand{\goto}{\rightarrow}
\newcommand{\eps}{\epsilon}

\usepackage{setspace}
\allowdisplaybreaks
\usepackage{xr}

\begin{document}

\title{\bf Optimal estimators and tests for reciprocal effects }

\author{Qunqiang Feng\thanks{Department of Statistics and Finance, University of Science and Technology of China, Hefei, 230000, China. \texttt{Email:} fengqq@ustc.edu.cn} \hspace{4mm} 
	\and Jiashun Jin\thanks{School of Statistics and Data Science, Southeast University, Nanjing, 211100, China. \texttt{Email:} jiashunmail@gmail.com, yrtian1997@gmail.com} \hspace{4mm} 
	\and Yaru Tian \footnotemark[2] \hspace{4mm} 
	\and Ting Yan \thanks{Department of Statistics, Central China Normal University, Wuhan, 430079, China. \texttt{Email:} tingyanty@mail.ccnu.edu.cn.} } 
\date{}

\maketitle
\let\thefootnote\relax\footnotetext{Author names are sorted alphabetically. Corresponding author: Jiashun Jin.}
\begin{abstract}
	\begin{spacing}{1.2}
		The $p_1$ model plays a  fundamental role in  modeling directed networks, where the reciprocal
		effect parameter $\rho$ is of special interest in practice. However,  due to nonlinear factors  in this model, how to estimate $\rho$ efficiently is a long-standing open problem.
		
		We tackle the problem by the cycle count approach.  The challenge is,   due to the nonlinear factors in the model, for any given type of generalized cycles, the expected count is a complicated
		function of many parameters in the  model, so it is unclear how to use cycle counts to estimate $\rho$.
		However, somewhat surprisingly, we discover that, among many types of generalized cycles with the same length,
		we can carefully pick a pair of them such that in the ratio between the expected cycle counts of the two types,
		the non-linear factors {\it cancel out nicely with each other}, and as a result, the ratio equals to $\mathrm{exp}(\rho)$ exactly.
		Therefore,  though the expected count  of cycles of any type is not tractable,
		the ratio between the expected cycle  counts of a (carefully chosen) pair of generalized cycles
		may have an utterly simple form.
		
		We study to what extent such pairs exist, and
		use our discovery to derive both an estimate for $\rho$ and a testing procedure for testing $\rho = \rho_0$. In a setting where we allow a wide range of reciprocal effects and a wide variety  of network sparsity and degree heterogeneity,
		we show that our estimator achieves the optimal rate and our test achieves the optimal phase transition. Technically,  first,   
		motivated by what we observe on real networks, we do not want to impose strong conditions on reciprocal effects, network sparsity,  and degree heterogeneity. Second,  our proposed statistic is a type of $U$-statistic,   the analysis of 
		which involves complex combinatorics and is error-prone. For these reasons, our analysis is long and delicate. 
	\end{spacing}

	\begin{spacing}{1.4}
	\end{spacing}
	
\end{abstract}

\section{Introduction}
\label{sec:intro}
Reciprocity is a fundamental principle in social psychology, where we believe that individuals tend to respond to the actions of others in a way that reflects the nature (positive/negative) of those actions, and people make choices based on what they can gain (profits,   trusts, feelings) from others in return.  Just like social duties,  reciprocity makes it possible to build sustainable and continuing relationships that bond different individuals in a community together.

In directed  networks, how to model  reciprocity is a problem of great interest.
The $p_1$ model by \citet{Holland}  is among  the most popular models  for directed networks. It is a powerful tool for characterizing reciprocation and has been widely applied in practice; see \cite{Goldenberg2010} for a survey.
Consider a directed network with $n$ nodes. Let $A$ be the
adjacency matrix, where $A_{ij}  = 1$ if  there is a (directed) edge
from  node $i$  to  node $j$.   As a convention, $A_{ii} = 0$. We  assume that the bivariate random variables $\{(A_{ij}, A_{ji}): 1 \leq i < j \leq n\}$ are independent (bivariate) Bernoulli, and that
for parameters $ \rho,\gamma,   \alpha_1, \beta_1, \ldots, \alpha_n, \beta_n$ and any $a, b \in \{0, 1\}$,
\begin{equation} \label{p1model1}
\mathbb{P}(A_{ij} = a, A_{ji} = b) =   K_{ij} \cdot \mathrm{exp}(a (\gamma + \alpha_i + \beta_j) + b (\gamma + \alpha_j + \beta_i) + a b \rho),  \;\;  1 \leq i \neq j \leq n,
\end{equation}
where  $K_{ij} = [1 + e^{\gamma + \alpha_i + \beta_j}  + e^{\gamma + \alpha_j + \beta_i} +
e^{2\gamma + \alpha_i + \beta_j + \alpha_j + \beta_i + \rho}]^{-1}$. For identifiability,
we assume
\begin{equation} \label{p1model2}
\alpha_1 + \alpha_2 + \ldots + \alpha_n =  \beta_1 + \beta_2 + \ldots + \beta_n = 0.
\end{equation}
Here, $\rho$ calibrates the strength of reciprocated edges, $\gamma$ calibrates the overall sparsity level,  and $\alpha_i$ and  $\beta_i$ calibrate the expansiveness and popularity of node $i$,  respectively.

The $p_1$ model can also be viewed as a special case of the Exponential Random Graph Model \cite{Goldenberg2010,  ERGM}.
It  includes the well-known $\beta$-model \cite{Chatterjee2011random} as a special case.
The $\beta$-model has been widely studied (e.g., \cite{Chen:2020, Yan:Leng:Zhu:2016, Yan2013clt}).
 Also, if we neglect the non-linear factors $K_{ij}$, then the $p_1$ model reduces to a special block model. The block models have received a lot of attention recently (e.g., \cite{CMM, fan2022simple,  SCORE-Review,  ke2022optimal, li2020network, Feng2015,   wang2017likelihood, zhao2012}). See Section  \ref{subsec:extension} and Section \ref{sec:Discu} for more discussion.

Note that real networks usually have severe degree heterogeneity \citep{SCORE-Review}  (meaning the degree, in-degree or out-degree,  of a node can be hundreds time larger than that of the other). Also,  the sparsity level of a network may range significantly from one occasion to another. For these reasons,
we allow $\gamma, \alpha_1,  \beta_1, \ldots, \alpha_n, \beta_n$ to take all possible values in their range of interest.

Let $\theta = (\rho, \gamma, \alpha_1, \beta_1,  \ldots,  \alpha_n, \beta_n)$. Due to the importance of reciprocity,
out of the $(2n+2)$ unknown parameters, $\rho$ (the reciprocity parameter) has received the most attention \cite{Goldenberg2010}. In many applications (e.g., link prediction,  citation prediction \cite{JBES}),
the parameter plays an important role.      Motivated by these, our primary interest of the paper is as follows.
\begin{itemize}
\setlength \itemsep{-.5 em}
\item {\it Optimal estimation of  $\rho$}.  How to estimate $\rho$ with the optimal rate of convergence.
\item {\it Optimal tests}.   Given a $\rho_0$, how to derive optimal tests for testing whether $\rho = \rho_0$.
\item {\it Optimal adaptivity}. The tests and estimators are optimal in a broad setting covering all possible levels of reciprocity and sparsity and allowing for severe degree heterogeneity.  \end{itemize}

Despite the popularity of the $p_1$ model in practice,  how to estimate $\rho$ and other parameters  is
a long-standing problem that has challenged us for over $40$ years \cite{Goldenberg2010, Rinaldo2013mle}.
The main reason is that, due to  factors $K_{ij}$ in Model (\ref{p1model1}),
the $p_1$ model is {\it nonlinear} (see Section \ref{sec:prelim}).
We now briefly review the literature on estimating $\rho$ (discussion for testing is similar and is omitted).

First, we may consider the Maximum Likelihood Estimation (MLE) approach.   However, as pointed out by \cite{Goldenberg2010}, MLE faces grand challenges:  ``A major problem with the $p_1$ and related models, ..., 
we have no consistency in results for the maximum likelihood estimates, and no simple way to test for $\rho=0$."
Also,  \citet[Theorem 1.5, supplement]{Rinaldo2013mle} proved that the MLE exists
if and only if the adjacency matrix $A$ in the $p_1$ model satisfies two hard-to-check conditions (in fact,  MLE and LRT do not exist in many cases; see Section \ref{sec:numeric}). Also, relatively, MLE is computationally slow and sensitive to outliers and model mis-specification.

Second, we may consider the spectral approach, which is popular in network analysis (e.g., \cite{SCORE-Review}).
However, spectral approaches are only appropriate for low-rank network models.
In a low-rank network model,  we can decompose the adjacency matrix $A$ as $A = \Omega + W$,
where the matrix $\Omega$ is low-rank and non-stochastic (containing the information of the parameters)
and the matrix $W$ contains the noise and secondary effect.
Unfortunately, due to the nonlinear factors $K_{ij}$ in Model (\ref{p1model1}),
the $p_1$ model is not a low-rank model (see Section \ref{sec:prelim} for more discussion), and it is unclear how to develop an efficient spectral approach to estimate $\rho$.

Last, we may consider the cycle count approach, which is also popular in network analysis. Consider
an undirected network with $K$ communities.  In recent works, \cite{gao2017testing, JinGC2018, JKL2021} developed a cycle count approach to  global testing where the goal is to test whether $K = 1$ or $K > 1$ and \cite{JKLW2022} developed a cycle count approach to estimating $K$.   Unfortunately, these works have been largely focused on (a) different problems, and (b) undirected networks with low-rank models, and it remains unclear how to extend their ideas to address our problem.
For these reasons, how to estimate $\rho$ is an interesting but largely unsolved problem.

\subsection{LCR: A novel cycle count approach to estimating $\rho$}
We wish to develop a cycle count approach for our problem.  The main challenge is,  for any $m \geq 3$, there are many different types of $m$-cycles (see below), and it is unclear which cycle count statistics
may contain the key information about $\rho$. 

To address the challenge, our idea is as follows. We first extend the notion of (regular) cycles for undirected networks to generalized cycles for
directed networks, with more complicated forms.   We then introduce  pseudo low-rank representation (PLR)  and  edge-encoding for the $p_1$ model.  Using PLR and edge encoding, for any type of length-$m$ generalized cycles, we can explicitly represent the expected number of such cycles with a big sum. The problem is,  however, due to the nonlinear factors $K_{ij}$ in (\ref{p1model1}),  it is hard to further
simplify the big sum  above
with a succinct formula  and relate it to our parameters.  Fortunately, we discover that, for some $m \geq 3$, we can construct   a pair of generalized $m$-cycles, Type-$a$ and Type-$b$, such that for an integer $c_0 \neq 0$ (known to us),
\begin{equation} \label{idea}
\frac{\mbox{Expected number of Type-$a$ generalized cycles}}{\mbox{Expected number of Type-$b$ generalized cycles}}  = e^{c_0 \rho}.
\end{equation}
Below, we show that there are $O(4^{2m})$ pairs of $(a,b)$, but out of them, only very few pairs satisfy (\ref{idea});  we need non-trivial efforts to find such pairs of $(a, b)$ (see Section \ref{sec:main}). 
Therefore, {\it  though it is hard to derive a simple formula for either the
numerator or denominator on the LHS of (\ref{idea}),  we can
derive a very simple formula for their ratio}.  This gives rise to
a convenient way to estimate $\rho$ and to test whether $\rho = \rho_0$ for a given $\rho_0$;
see Section \ref{sec:main}.

As far as we know,  the discovery in (\ref{idea})  is new.
We have the following contributions. (1).  When $m$ is even, we can always construct a pair of length-$m$ generalized cycles such that (\ref{idea}) is satisfied. However, when $m$ is odd, such a pair does not exist. (2).  Using the above discovery, we  propose the logarithmic cycle-count ratio (LCR) as a novel approach to estimating $\rho$. In a broad setting where $\rho$ may vary in its whole range of interest, and which covers all interesting sparsity levels and allows severe degree heterogeneity, we show that LCR is asymptotical minimax (and thus is optimal).
(3).  We propose a new way to estimate the variances of LCR estimators, and
use the results to develop the LCR tests for testing $\rho = \rho_0$.  In the same broad setting,  we show that the LCR  tests  have $N(0,1)$ as the limiting null  and  achieve  the optimal phase transition (similar to the minimax framework, phase transition is a theoretical framework for optimality \cite{DJ15, phase}). In practice,  an explicit limiting null is especially useful (e.g., for computing $p$-values). 
In summary, our approach is  new, and for the first time, we derive  optimal tests and estimators  for the reciprocal parameter $\rho$ in a broad setting (by far,  consistent estimators
for $\rho$ in a broad setting are not available).

{\bf Organizations and notations}.  Section \ref{sec:prelim} presents  the Pseudo Low-rank Representation (PLR) for the $p_1$ model.
Section \ref{sec:main} contains our main results.   Section \ref{sec:numeric} contains numerical studies.   Section \ref{sec:proof}
 sketches the proof ideas.  Section \ref{sec:Discu} is a short discussion.  We use $C > 0$ to denote a generic constant which may vary from one occasion to another.
For a vector $x=(x_1,\ldots,x_n)' \in\mathbb{R}^n$,   $\|x\|_q$ denotes the
$\ell^q$-norm of $x$, and $x_{max}$ and $x_{min}$ denote the largest and smallest entry of $x$, respectively.
For any positive sequences $\{a_n\}_{n  = 1}^{\infty}$ and $\{b_n\}_{n = 1}^{\infty}$,
we say $a_n \sim b_n$ if $\lim_{n \goto \infty} (a_n / b_n) = 1$ and 
$a_n \asymp b_n$ if there are constants $c_1 > c_0 > 0$ such that
$c_0 \leq (a_n / b_n) \leq c_1$ for sufficiently large $n$.    For $A, B \in \mathbb{R}^{n,m}$ ($m \geq 1$),  $A \circ B \in \mathbb{R}^{n,m}$ denotes the matrix of  entry-wise product. For any $x, y \in \mathbb{R}^n$,   $(x, y)$ denotes their inner product.  

\section{Preliminary:  The PLR and generalized cycles} \label{sec:prelim}
We introduce the Pseudo Low-rank Representation (PLR) and edge coding for   model (\ref{p1model1}).  Let $\circ$ be the Hadamard or entry-wise product \cite{HornJohnson}.  Let  $K$ be the $n \times n$ matrix $K = (K_{ij})$.
Let $\mu, \nu$ be the $n \times 1$ vectors satisfying  $\mu_i = e^{\gamma/2 + \alpha_i}$ and $\nu_i = e^{\gamma/2 + \beta_i}$, $1 \leq i \leq n$, and let $\eta = \mu \circ \nu$.  Let $\widetilde{\Omega} = \mu \nu' + e^{\rho} \eta \eta'$ and $\Omega = K \circ \widetilde{\Omega}$.  By (\ref{p1model1}), it is seen that  $\Omega_{ij} = \mathbb{P}(A_{ij} = 1)$, $1 \leq i \neq j \leq n$. Let $\diag(\Omega)$ be the $n \times n$ diagonal matrix where the $i$-th diagonal entry is $\Omega_{ii}$. Let $W$ be the $n
\times n$ matrix where $W_{ij} = A_{ij} - \Omega_{ij}$ if $i \neq j$ and $W_{ij} = 0$ otherwise.
It follows  that
\begin{equation} \label{plowrankmodel1}
A = \Omega - \diag(\Omega) + W, \qquad \mbox{where $\Omega =  K \circ \widetilde{\Omega}$ with $\widetilde{\Omega} = [\mu  \nu' + e^{\rho} \eta \eta']$}.
\end{equation}
Let $J_n \in \mathbb{R}^{n,n}$ be the matrix of all ones and $\tilde{J}_n = J_n - I_n$, where $I_n$ is the identity matrix.
Introduce four $n \times n$ matrix as follows (assuming $i \neq j$ in the paragraph below):  (1).  $A^{11} = A \circ A'$. Note that $A_{ij}^{11} = 1$ if and only if there is a double edge between $i$ and $j$; (2).  $A^{00} = (\tilde{J}_n - A) \circ (\tilde{J}_n - A')$ ($A_{ij}^{00} = 1$ if and
only if there is no edge between $i$ and $j$); (3).  $A^{10} = A \circ (\tilde{J}_n - A')$ and its transpose $A^{01} = (A^{10})'$. Note that $A_{ij}^{10} = 1$ if and only if there is a (directed) edge from $i$ to  $j$ but there is no edge from $j$ to $i$; similar for $A^{01}$. Note also that $A^{11}$ and $A^{00}$ are symmetrical but $A^{10}$ and $A^{01}$ may not.    Let $\widetilde{\Omega}^{00} = J_n$,  $\widetilde{\Omega}^{10} = (\widetilde{\Omega}^{01})' = \mu \nu'$ and $\widetilde{\Omega}^{11} = e^{\rho}\eta \eta'$. For $a \in \{00, 11, 10, 01\}$,
let $\Omega^a = K \circ \widetilde{\Omega}^a$. Also, similarly,  let
$\diag(\Omega^a)$ be the $n \times n$ diagonal matrix where the $i$-th diagonal entry is $\Omega_{ii}^a$, and let  $W^a$ be the $n \times n$ matrix where $W_{ij}^a = A_{ij}^a - \Omega_{ij}^a$ if $i \neq j$ and $W^a_{ij} = 0$ otherwise. It follows that
\begin{equation} \label{plowrankmodel3}
A^a = \Omega^a  - \diag(\Omega^a) + W^{a}  \;\;\;  \mbox{and}  \;\;\;   \Omega^a = K \circ \widetilde{\Omega}^a, \qquad \mbox{for any $a \in \{00, 11, 10, 01\}$}.
\end{equation}
We call (\ref{plowrankmodel1})-(\ref{plowrankmodel3}) the {\it pseudo low-rank representation (PLR)} of the $p_1$ model, for if the nonlinear matrix $K$ is dropped, then matrices $\Omega, \Omega^{00}, \Omega^{10}, \Omega^{01}, \Omega^{11}$ are all low-rank. 
 
We now introduce generalized cycles and edge coding.  Fix $m \geq 3$.  A generalized $m$-cycle
contains $m$ distinct nodes, $i_1, i_2, \ldots, i_m$, and many edges. 
For any pair of nodes $i, j$ in $\{i_1, i_2, \ldots, i_m\}$,   there are $4$ different types of edges: no edge, a directed edge from $i$ to $j$, a directed edge from
$j$ to $i$, and a double edge;  we encode them by $\{00, 10, 01, 11\}$.  
Therefore, to consider a specific type of length-$m$ generalized cycle, we need to specify
the edge pattern first. Introduce
$S_m = \{a = (a_1, a_2, \ldots, a_m):  \mbox{where each $a_i$ takes values from $\{00, 10, 01, 11\}$}\}$. 
\begin{definition}
For $m \geq 3$, a generalized $m$-cycle with edge pattern prescribed by an $a \in S_m$
consists of $m$ distinct nodes $i_1, i_2, \ldots, i_m$,  where the edge type
between $i_1$ and $i_2$,  $i_2$ and $i_3$, $...$,  and $i_m$ and $i_1$ are
$a_1, a_2, ..., a_m$, respectively (we call it an $m$-cycle of Type-$a$).
\end{definition}
For example, when $m = 4$, we have $4^4 = 256$ different generalized $4$-cycles (quadrilaterals) and $42$ different non-isomorphic generalized quadrilaterals.  See Tables A1 of the supplement.

By the PLR and edge-coding, we can represent the number of any given type of
generalized cycles succinctly by a big sum.   Fix $m \geq 3$ and $a \in S_m$. By the PLR, the number of generalized $m$-cycles of Type-$a$ is $\frac{1}{C_m(a)} \sum_{i_1, i_2, \ldots, i_m (dist)}  A_{i_1i_2}^{a_1}  A_{i_2i_3}^{a_2}  \ldots A_{i_m i_1}^{a_m}$.
This is because the value of $A_{i_1i_2}^{a_1}  A_{i_2i_3}^{a_2}  \ldots A_{i_m i_1}^{a_m}$ is either $0$ and $1$,
and it is  $1$ if and only if there is a  type-$a_1$ edge from $i_1$ to $i_2$,
a type-$a_2$ edge from $i_2$ to $i_3$, $\ldots$, and a type-$a_m$ edge from $i_m$ to $i_1$.
Here, $C_m(a)$ is a combinatoric number representing the number of times we have counted each cycle repeatedly.
This motivates us to consider the statistic (where $G$ stands for generalized cycle)
\begin{equation} \label{idea1}
G_{n, m}(a) = \sum_{i_1, i_2, \ldots, i_m (dist)}  A_{i_1i_2}^{a_1}  A_{i_2i_3}^{a_2}  \ldots A_{i_m i_1}^{a_m}.
\end{equation}
Since $i_1, i_2, \ldots, i_m$ are distinct,  the $m$ random variables $A_{i_1i_2}^{a_1},  A_{i_2i_3}^{a_2}, \ldots, A_{i_m i_1}^{a_m}$ are independent of each other.
Therefore, if we let $g_{n, m}(a) = \mathbb{E}[G_{n, m}(a)]$, then it follows that
\begin{equation} \label{idea2}
g_{n, m}(a) =  \mathbb{E}[G_{n, m}(a)] =   \sum_{i_1, i_2, \ldots,  i_m (dist)}    \Omega_{i_1i_2}^{a_1}  \Omega_{i_2i_3}^{a_2}  \ldots \Omega_{i_m i_1}^{a_m}.
\end{equation}
Below we use these big sums to introduce our main ideas and results.

\section{Main results} \label{sec:main}
In this section, we propose an optimal approach to estimating $\rho$ and to testing
whether $\rho = \rho_0$. 
\subsection{A novel approach for estimating $\rho$ using a cancellation trick} \label{subsec:trick}
Fix an $m \geq 3$ and $a \in S_m$ and let $G_{n,m}(a)$ and $g_{n,m}(a)$ be as in (\ref{idea1})-(\ref{idea2}).
To construct an estimator for $\rho$, as the first attempt, we hope  to find $a \in S_m$ such that
$\rho$ is a simple function of $g_{n, m}(a)$. Note that we expect to see that  $G_{n, m}(a) \approx g_{n,m}(a)$.  So if we can express $\rho$ as a simple function of $g_{n, m}(a)$, then a simple estimator for $\rho$ arises. Unfortunately, due to the nonlinear factors $K_{ij}$,   such an approach does not work.
In fact, by (\ref{plowrankmodel1}) and (\ref{idea2}),
\[
g_{n,m}(a) =  \sum_{i_1, i_2, \ldots,  i_m (dist)}    \Omega_{i_1i_2}^{a_1}  \Omega_{i_2i_3}^{a_2}  \ldots \Omega_{i_m i_1}^{a_m} =  \sum_{i_1, i_2, \ldots,  i_m (dist)}  K_{i_1 i_2}  \ldots K_{i_m i_1}   \widetilde{\Omega}_{i_1i_2}^{a_1}    \ldots \widetilde{\Omega}_{i_m i_1}^{a_m}.
\]
Due to the nonlinear factors $K_{i_1i_2} \ldots K_{i_m i_1}$ in each term in the big sum, it is hard to further simplify the RHS and derive an explicit formula in terms of $\theta = (\rho, \gamma, \alpha_1, \beta_1, \ldots,
\alpha_n, \beta_n)$.

To overcome the difficulty, our idea is as follows.  Fix $a, b \in S_m$ and consider the ratio $G_{n, m}(a) / G_{n, m}(b)$.
Similarly, we expect to see $G_{n,m}(a) / G_{n, m}(b) \approx g_{n,m}(a) / g_{n,m}(b)$,
where
\begin{equation} \label{idea3}
\frac{g_{n,m}(a)}{g_{n,m}(b)} = \frac{\sum_{i_1, i_2, \ldots,  i_m (dist)}  K_{i_1 i_2}  K_{i_2 i_3} \ldots K_{i_m i_1}   \widetilde{\Omega}_{i_1i_2}^{a_1}  \widetilde{\Omega}_{i_2 i_3}^{a_2}   \ldots \widetilde{\Omega}_{i_m i_1}^{a_m}}{\sum_{i_1, i_2, \ldots,  i_m (dist)}  K_{i_1 i_2} K_{i_2 i_3}  \ldots K_{i_m i_1}   \widetilde{\Omega}_{i_1i_2}^{b_1}  \widetilde{\Omega}_{i_2 i_3}^{b_2}  \ldots \widetilde{\Omega}_{i_m i_1}^{b_m}}.
\end{equation}
Note that (a)  the non-linear factors $K_{i_1i_2} K_{i_2 i_3} \ldots K_{i_m i_1}$ {\it do not depend on}  $a$ and $b$ and are the same in the numerator and denominator, and (b)  the matrices  $\widetilde{\Omega}^{a_i}$ and $\widetilde{\Omega}^{b_i}$ are relatively simple and  easy-to-handle. Therefore, if we can {\it find a pair of $(a, b)$ where we can cancel the nonlinear factors $K_{i_1i_2} K_{i_2 i_3} \ldots K_{i_m i_1}$ in each term of the numerator and the denominator}, then we are able to derive a simple and closed-form formula for $g_{n,m}(a) / g_{n,m}(b)$.

Somewhat surprisingly, this turns out to be possible. In detail,  if we pick $a, b \in S_m$ such that for an explicit and known number $c_0(a, b) > 0$,
\begin{equation} \label{critical}
\widetilde{\Omega}_{i_1 i_2}^{a_1}  \widetilde{\Omega}_{i_2 i_3}^{a_2}  \ldots \widetilde{\Omega}_{i_m i_1}^{a_m} =  e^{c_0(a, b) \rho}  \cdot \widetilde{\Omega}_{i_1 i_2}^{b_1}  \widetilde{\Omega}_{i_2 i_3}^{b_2}  \ldots \Omega_{i_m i_1}^{b_m}, \qquad \mbox{for all distinct $i_1, i_2, \ldots, i_m$},
\end{equation}
then each nonlinear factor $K_{i_1i_2} K_{i_2 i_3} \ldots K_{i_m i_1}$ in  the big sum of the numerator of (\ref{idea3}) cancels with the same factor in the denominator,      and so $g_{n,m}(a) / g_{n,m}(b) = e^{c_0(a, b) \rho}$; see Lemma \ref{lemma:ab}.

\begin{figure}[htb!]
\centering
\includegraphics[height = 2.2 in]{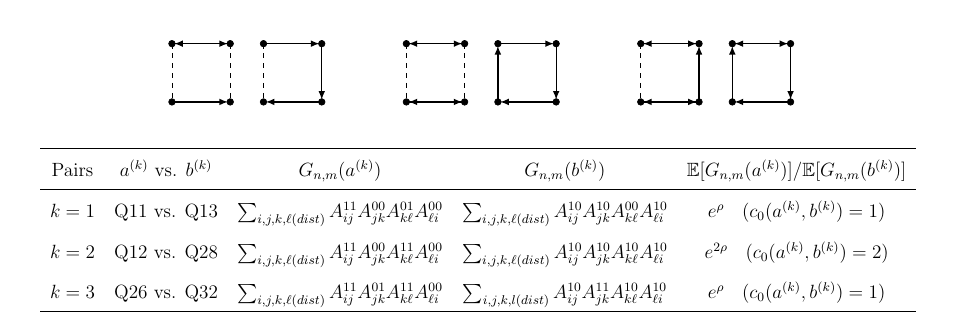} 
\caption{The three non-isomorphic pairs of $(a,b)$  that satisfy condition (\ref{critical}) for $m = 4$ (e.g.,  Q11 is the
11-th  generalized $4$-cycle (quadrilateral) in Table A1 of the supplement). 
Top (from left to right): generalized quadrilaterals for $(Q11, Q13)$, $(Q12, Q28)$, and $(Q26, Q32)$.}
\label{tab:pairs}
\end{figure}

It remains to construct $a, b \in S_m$ such that (\ref{critical}) holds.
When $m$ is odd, such a construction does not exist. However,
for an even number $m \geq 4$, such a construction always exist.
 Recall that  $a_i, b_i$ take values in
$\{00, 01, 10, 11\}$.    Lemma \ref{lemma:ab} is proved in the supplement.
\begin{lemma} \label{lemma:ab}
If (\ref{critical}) hold,  then $g_{n,m}(a) / g_{n, m}(b) =    e^{c_0(a, b)  \rho}$.
When $m \geq 3$ and is odd, no pair of $(a, b)$ satisfies condition (\ref{critical}).
When $m \geq 3$ and is even, such pairs exist if we let $N = m/2$,  $G_{n, m}(a) = \sum_{i_1, \ldots, i_m (dist)} A_{i_1 i_2}^{1x_1} A_{i_2 i_3}^{y_10} \ldots A_{i_{m-1} i_m}^{1x_N} A_{i_m i_1}^{y_N0}$,   and
$G_{n,m}(b) = \sum_{i_1, \ldots, i_m (dist)} A_{i_1 i_2}^{0x_1} A_{i_2 i_3}^{y_11} \ldots A_{i_{m-1} i_m}^{0x_N} A_{i_m i_1}^{y_N1}$, 
where $x_1, y_1, \ldots, x_N, y_N \in \{0,1\}$ satisfy $c_0(a, b):  = (x_1   + \ldots +  x_N) - (y_1 +   \ldots +  y_N) > 0$. 
Moreover, this gives all possible pairs of $(a, b)$:   for any  pair  $(a', b')$ satisfying (\ref{critical}), there is a pair of $(a, b)$ given above
 such that the generalized cycles associated with
$a'$ and $b'$ are isomorphic to those associated with $a$ and $b$, respectively.
\end{lemma}
When $m = 4$, there are a total $3$ non-isomorphic pairs of $(a, b)$ satisfying (\ref{critical}); see Figure  \ref{tab:pairs}. For each $1 \leq k \leq 3$, there is a pair of $(a,b)$ given by Lemma \ref{lemma:ab} with $m =4$ such that the generalized cycles corresponding to $a^{(k)}$ and $b^{(k)}$ are isomorphic to the
generalized cycles corresponding to $a$ and $b$, respectively.  For example, the generalized cycles defined by $(a^{(1)}, b^{(1)})$ in Figure \ref{tab:pairs}  are isomorphic to $A^{1x_1} A^{y_10} A^{1x_2} A^{y_20}$ and
$A^{0x_1} A^{y_11} A^{0x_2} A^{y_21}$,   respectively, with $(x_1,  y_1,x_2, y_2) = (1, 0, 0, 0)$. See also 
Figure \ref{tab:pairs} for the geometrical shapes of these cycles. 
 
We now discuss how to estimate $\rho$.  By Remark 2 below, we focus on $m = 4$. 
When $m = 4$,  using any of the $3$ pairs of $(a, b)$ in Figure \ref{tab:pairs},     we can successfully cancel the nuisance nonlinear factors
$K_{i_1i_2} K_{i_2 i_3} \ldots K_{i_m i_1}$, and a  simple estimator for $\rho$ arises from
\[
G_{n,m}(a) / G_{n,m}(b)   \approx  g_{n,m}(a) / g_{n,m}(b)  = \mathrm{exp}(c_0(a, b)  \rho),  \qquad (\mbox{note: $m = 4$ and  $c_0(a, b)$ is known}). 
\]
This gives rises to three different estimators for $\rho$. In Section \ref{subsec:otherpairs},  we show that among the  three estimators,  the first one has the fastest rate.  For this reason,  we focus   on the first pair of $(a, b)$  below. 
Write for short $(a, b) = (a^{(1)}, b^{(1)})$.  Introduce ($Q$ stands for quadrilateral)
\[
Q_n(a)  \equiv G_{n, 4}(a) =  \sum_{i, j, k,  \ell (dist)} A_{ij}^{11} A_{jk}^{00}   A_{k \ell}^{01}  A_{\ell i}^{00}, \;\;\;   Q_n(b) \equiv G_{n, 4}(b) =  \sum_{i, j, k,  \ell (dist)} A_{ij}^{10} A_{jk}^{10}   A_{k \ell}^{00}  A_{\ell i}^{10}.
\]
Let $H(\cdot; t)$ be the hard-thresholding operator where $H(X; t) = X$ if $|X| \leq t$ and $H(X; t) = \mathrm{sign}(X) \cdot t$ otherwise.  We propose the {\it Logarithmic Cycle-count Ratio (LCR)} estimate by
\begin{equation} \label{DefineLCR}
\hat{\rho}^* = H(\hat{\rho}; 2 \log(n)), \qquad \mbox{with} \; \hat{\rho} = \log(Q_n(a) / Q_n(b)).
\end{equation}
The hard-thresholding makes it easier to control rare events, but practically it does not  make any difference:  the range of interest  is  $n^{-2} \ll e^{\rho} \ll n^2$ (and so $|\rho| < 2 \log(n)$); see Section \ref{subsec:condition}.

Consider now the computing cost of $\hat{\rho}^*$. By the PLR, when $\{i, j, k, \ell\}$ are not distinct,
$A_{ij}^{11} A_{jk}^{00}  A_{k\ell}^{01}  A_{\ell i}^{00}  =    A_{ij}^{10}  A_{jk}^{10}  A_{k\ell}^{00} A_{\ell i}^{10} = 0$.
Using this, we have Lemma \ref{lemma:trace} (the proof  is skipped).  
\begin{lemma}  \label{lemma:trace}
$Q_n(a) = \mathrm{trace}(A^{11} A^{00} A^{01} A^{00})$ and $Q_n(b) = \mathrm{trace}(A^{10} A^{10} A^{00} A^{10})$.
\end{lemma}
The main computing cost of $\hat{\rho}^*$ comes from computing
$Q_n(a)$ and $Q_n(b)$, but the complexity of computing $Q_n(a)$ and $Q_n(b)$
are $O(n^2 d)$, where $d = \max_{1 \leq i \leq n} d_i$, and $d_i$ is the maximum
of in-degree, out-degree, and double-edge degree of node $i$.
For very sparse networks, $d$ is small, and the complexity is
only larger than $n^2$ by a comparably small factor.

{\bf Remark 1}.   Our paper is connected to some works in network analysis (e.g., \cite{gao2017testing, JinGC2018, JKL2021, JKLW2022, CCDS}), which also use cycle count statistics, but there are  major differences. First, our paper is for directed networks, where there are many types of $m$-cycles (and our statistic is the ratio of two types of $m$-cycles).    
The above cited works are for undirected networks, where we only have one type of $m$-cycles. 
Second, our cycle count statistic is applied to the original adjacency matrix $A$, 
while those in \cite{JKL2021,JKLW2022} are applied to a centerdized adjacency matrix $A - \widehat{\Omega}$.   
 Last but not the least, our paper is for estimating the reciprocal effects, while \cite{gao2017testing, JinGC2018, JKL2021, JKLW2022} focused on the testing and estimation problems  associated with the number of network communities.   For these reasons, our idea (including the cancellation trick,  the construction of $(a, b)$, and the proposed approaches below) is new and was not discovered before.

{\bf Remark 2} (The case of larger $m$). By Lemma \ref{lemma:ab},  for any even number  $m \geq 4$,
we can construct a pair of $(a,b) \in S_m$ such that (\ref{critical}) holds. This means we can
have many similar estimators as $\hat{\rho}^*$. However, the estimator for a larger $m$ is {\it much more
tedious} to analyze. Since $\hat{\rho}^*$ is already rate optimal (see below),   we leave further investigation to the future.

{\bf Remark 3}.
Our approach is a {\it non-parametric direct approach}  for it does not require estimating the other $(2n+1)$ parameters in the   model.
A direct approach is especially desirable. 
First,  it is unclear what is the best estimates for   other  parameters.   Second, if we first estimate  other  parameters and then use the results to estimate $\rho$,   then we may have non-optimal results, for
the rates of  estimating other  parameters may be slower than that of estimating $\rho$.
Last,  many times, a direct approach is faster and conceptually simpler. See Section \ref{subsec:extension}.  

{\bf Remark 4}.  The LCR  has a closed-form and is optimal in a broad setting (see below).
In comparison, the MLE does not have a closed-form and is analytically challenging:  it remains unclear how to check whether the MLE exists for a given data set and 
how to analyze the variance of the MLE  (e.g., \cite{Goldenberg2010, Rinaldo2013mle}). In fact, numerically, we frequently observe that the MLE does not exist,  especially when the networks are sparse, but theoretically, we have limited understanding on whether the MLE exists or not (e.g., \cite{Rinaldo2013mle}).
Also,  recall that the complexity of LCR is $O(d n^2)$, where $d$ is the maximum degree.  The complexity of MLE depends on which algorithm we use for implementation. If we use that by \cite{Holland:Leinhardt:1981}, then the complexity is $O(rn^2)$, where $r$ is the number of iterations for reaching the predefined error value.  Our study suggests that $r$ is more sensitive to degree heterogeneity than sparsity,  and can be large for networks with severe degree heterogeneity. Therefore, for sparse networks with severe 
degree heterogeneity, the complexity of LCR  is smaller than that of MLE. For dense networks with moderately degree heterogeneity, we have the opposite.   See Section \ref{sec:numeric} for more comparison.



\subsection{Regularity conditions} \label{subsec:condition}
In this paper, we let $n \goto \infty$, and allow $\theta = (\rho, \gamma, \alpha_1,\beta_1,\ldots, \alpha_n, \beta_n)$
to depend on $n$ as $n$ varies (but for notational simplicity, we still write $\rho = \rho_n, \gamma= \gamma_n,
\alpha_i = \alpha_{in}, \beta_i = \beta_{in}$, without the subscripts).
By the PLR in  Section \ref{sec:prelim}, for $1 \leq i \leq n$,
$\mu_i = e^{\alpha_i + \gamma/2}$, $\nu_i = e^{\beta_i + \gamma/2}$, and $\eta_i = \mu_i \nu_i$.
Let $\mu_{max}$ and $\mu_{min}$ denote the largest and smallest entry in $(\mu_1, \ldots, \mu_n)$ respectively (similar for
$\nu_{max}$ and $\nu_{min}$). For any $i \neq j$, the probability that there is no edge between them is $K_{ij}^{-1} = [1 + \mu_i \nu_j  + \mu_j \nu_i  + e^{\rho} \eta_i \eta_j]^{-1}$.  Since most real networks are sparse, we assume
\begin{equation} \label{condition1}
(\mbox{G1}): \qquad \mu_{max} \goto 0, \qquad  \nu_{max} \goto 0,  \qquad e^{\rho/2} \eta_{max} \goto 0,
\end{equation}
so for each node, the in-degree, out-degree, and degree of double edges are $o(n)$.\footnote{The variance of a cycle count statistic usually has many terms, with a cumbersome form. 
Fortunately, it was discovered (e.g., \cite{JKL2021,JKLW2022}) that with a mild sparsity condition, then asymptotically,  the variance reduces to a simple form.   G1 is part of such a condition: It is largely for technical reasons, and can be relaxed.\\} 
Now, when (\ref{condition1}) holds,  the total degree of double edges
is $\sim  e^{\rho} \|\eta\|_1^2$, so it is necessary to assume $e^{\rho/2} \|\eta\|_1 \goto \infty$.
Also,  for $1 \leq i \leq n$, the in-degree and out-degree of node $i$
are $\sim \nu_i \|\mu\|_1$ and $\sim \mu_i \|\nu\|_1$, respectively. Recall that $\eta_i = \mu_i \nu_i$. It is therefore    necessary to assume $\|\eta\|_1 \goto \infty$.  In light of this, we assume
\begin{equation} \label{condition2}
(\mbox{G2}): \qquad e^{\rho/2} \|\eta\|_1 \goto \infty, \qquad \|\eta\|_1 \goto \infty.
\end{equation}
By explanations above, both items in (\ref{condition2}) are mild conditions.  

(G1) and (G2) are the general conditions we need in the paper.  In some cases, we also need two {\it specific}  conditions. 
First, note that that both $Q_n(a)$ and $Q_n(b)$ are  big sums of many terms, where in order 
to accommodate severe degree heterogeneity, different terms may have very different magnitudes. 
In order for the central limit theorem (CLT)  and well-known large-deviation inequalities (e.g., \cite{Durrett2010Probability})  to work, 
we need a mild condition: 
\begin{equation} \label{condition3}
(\mbox{SP1}): \qquad  (\mu_{max}^2\nu_{max}^2)  / [(\mu, \eta)(\nu, \eta)]  \goto 0,   
\end{equation}
where $(x, y)$ denote the inner product of $x$ and $y$.  For CLT to work in $Q_n(a)$ and $Q_n(b)$, 
(SP1) is the minimum condition we need and 
can not be further relaxed.  The condition is also only mild (e.g., in the MDH case in Remark 5 below, 
$(\mu_{max}^2 \nu_{max}^2)  / [(\mu, \eta)(\nu, \eta)]  \asymp (1/[\|\mu\|_1 \|\nu\|_1])$, so (SP1) is implied by (G2)). 
Second,  the range of interest for $e^{\rho}$ is $n^{-2} \ll e^{\rho} \ll  n^2$ and it can be very large in some cases.  
Introduce the quantity $\tilde{\rho} = \tilde{\rho}(\mu,  \nu, \eta) = \log(\|\eta\|_1^2 / [(\mu, \eta)(\nu, \eta)])$. Note that by (\ref{condition1})-(\ref{condition2}), $\tilde{\rho} \geq -\log(\mu_{max} \nu_{max}) \goto \infty$.
We discuss two cases for $\rho$.
\begin{itemize}
\setlength \itemsep{-.5 em}
\item Case (S). In this case, $e^{\rho} \ll  e^{\tilde{\rho}}$. In this case, $e^{\rho}$ is small or moderately large at most (note: it includes the case of $\rho \leq C$).  Practically, this is  the most interesting case.
\item Case (L). In this case, $e^{\rho} \geq C e^{\tilde{\rho}}$ and $e^{\rho}$  is relatively large.
Sometimes,  we may further divide it to (L1) and (L2), where $e^{\rho}  \asymp e^{\tilde{\rho}}$ and
$e^{\rho} \gg e^{\tilde{\rho}}$, respectively.
\end{itemize}
The analysis for Case (L) is more delicate and requires a specific condition below: 
\begin{equation} \label{condition4}
(\mbox{SP2}): \qquad \|\mu\|_1 \asymp \|\nu\|_1,   \qquad \mbox{in (and only in)  Case (L)}.
\end{equation}
(SP2) requires that the total in- and out-degrees are at the same order. This is a mild condition: it always holds in the MDH  case (see Remark 5)   and is only needed
for Case (L).  

Conditions (\ref{condition1})-(\ref{condition4}) are all the main regularity conditions  in the paper,  except for
that in some places, we may slightly strengthen (\ref{condition2}) with a multi-$\log(n)$ factor;  see (\ref{conditionstr1}) below.  Here,  (\ref{condition1})-(\ref{condition2}) are general conditions which we use in most places, and (\ref{condition3})-(\ref{condition4}) are  specific conditions which we only use in some places. 
These conditions define a broad class of settings where we allow (a) $e^{\rho}$ to fully range in $(n^{-2}, n^2)$, (b) a wide range of sparsity, where
the total degrees may fully vary in $[1, n]$ (up to a  multi-$\log(n)$ factor  on both ends), and (c) severe degree heterogeneity:
$\mu_i$, $\nu_i$, and $\eta_i$
may be at very different magnitudes for different $1 \leq i \leq n$.

{\bf Remark 5}. (The moderate degree heterogeneous (MDH) case).   In most of our results,  we allow severe degree heterogeneity (SDH).   However, for illustration purpose,  we may also consider the MDH case where there is a constant $c_0 > 0$ such that
\begin{equation} \label{condition:moderate}
\max_{\{1 \leq i \leq n\}} \{  \max \{ |\alpha_i|,  |\beta_i|\} \} \leq c_0.
\end{equation}
In such a case, $\max\{\frac{\mu_{max}}{\mu_{min}}, \frac{\nu_{max}}{\nu_{min}}, \frac{\eta_{max}}{\eta_{min}} \} \leq C$ and condition (\ref{condition4}) holds. 
The average  in- or out- degree is $O(n e^{\gamma})$ and the total $\#$ of double edges is $O(n e^{\rho + 2 \gamma})$. Conditions (\ref{condition1})-(\ref{condition3}) reduce to $1 \ll ne^{\gamma} \ll n$  and $1 \ll  n^2 e^{\rho + 2\gamma}  \ll n^2$. This covers all sparsity levels in the range of interest.

\subsection{The Signal-to-Noise Ratio (SNR) of the LCR estimate $\hat{\rho}^*$} \label{subsec:SNR}
Recall that $\hat{\rho}^* = H(\hat{\rho}; 2 \log(n))$ where $\hat{\rho} = \log(Q_n(a) / Q_n(b))$.
  Introduce
\begin{equation} \label{DefineUn}
U_n(\rho)  =  Q_n(a) -  e^{\rho}Q_n(b).
\end{equation}
Recall that $\mathbb{E}[Q_n(a)] = e^{\rho} \mathbb{E}[Q_n(b)]$. By basic calculus, we expect to see
$\hat{\rho} - \rho = \log(\frac{Q_n(a)}{e^{\rho} Q_n(b)})    =   \log(1 + \frac{U_n(\rho)}{e^{\rho} Q_n(b)}) \approx
\frac{U_n(\rho)}{e^{\rho} \mathbb{E}[Q_n(b)]}
= \frac{U_n(\rho)}{\mathbb{E}[Q_n(a)]}$. Therefore, the key to analyze $\hat{\rho}^*$ is to analyze   
the Signal-to-Noise Ratio (SNR), defined  by
\begin{equation} \label{DefineSNR}
\mathrm{SNR} =  \mathbb{E}[Q_n(a)] / [\mathrm{Var}(U_n(\rho))]^{1/2}. 
\end{equation}
We now analyze the SNR. 
For $1 \leq i \neq j \leq n$, let $r_{ij} = \sum_{k, \ell \notin \{i, j\}, k \neq \ell} \Omega_{jk}^{00} \Omega_{k\ell}^{01} \Omega_{\ell i}^{00}$,  $s_{ij} = \sum_{k, \ell \notin \{i, j\}, k \neq \ell} \Omega_{jk}^{00} \Omega_{k\ell}^{11} \Omega_{\ell i}^{00}$, and
$t_{ij} = \sum_{k, \ell \notin \{i, j\}, k \neq \ell} (\Omega_{jk}^{10} \Omega_{k\ell}^{10} \Omega_{\ell i}^{00} + \Omega_{jk}^{00} \Omega_{k\ell}^{10} \Omega_{\ell i}^{10} + \Omega_{jk}^{10} \Omega_{k \ell}^{00} \Omega_{\ell i}^{10})$. Recall that
  $\theta = (\rho, \gamma, \alpha_1, \beta_1, \ldots, \alpha_n, \beta_n)$.
  Introduce
\begin{equation} \label{DefineV}
V_n(\rho) = \sum_{i \neq j} [2 r_{ij}^2 \Omega_{ij}^{11} + (s_{ij}  - e^{\rho} t_{ij})^2 \Omega_{ij}^{10}],  \;\;\;    r_n(\theta)  =  \max\bigl\{1/[e^{\rho} \|\eta\|_1^2], \; (\mu, \eta)(\nu, \eta) /\|\eta\|_1^4\bigr\}.
\end{equation}
\begin{theorem} \label{thm:snr1}
Consider the LCR estimator $\hat{\rho}^*$ in (\ref{DefineLCR}) and suppose conditions (\ref{condition1}), (\ref{condition2}), (\ref{condition4}) hold,
where $(\mu, \nu, \eta)$ are defined in (\ref{plowrankmodel1})  as in the PLR of the $p_1$ model.
As $n \goto \infty$,
\begin{itemize}
\item $\mathbb{E}[Q_n(a)]  = e^{\rho} \mathbb{E}[Q_n(b)]  = (1 + o(1))  \cdot  e^{\rho} \|\eta\|_1^2 \|\mu\|_1 \|\nu\|_1$.
\item $\mathrm{Var}(U_n(\rho))  = (1+o(1))V_n(\rho)
 =   (1 + o(1)) [2e^{\rho}\|\eta\|_1^2\|\mu\|_1^2\|\nu\|_1^2 +e^{2\rho}(\mu,\eta)(\nu,\eta)\|\mu\|_1^2\|\nu\|_1^2 +3e^{2\rho}(\mu,\eta)\|\eta\|_1^2\|\mu\|_1\|\nu\|_1^2  +3e^{2\rho}(\nu,\eta)\|\eta\|_1^2\|\mu\|_1^2\|\nu\|_1 -3e^{2\rho}\|\eta\|_1^4\|\mu\|_1\|\nu\|_1]$.
\item Moreover,
$(1/2 + o(1)) \leq \mathrm{Var}(U_n(\rho))   / [\mathrm{Var}(Q_n(a)) + e^{2\rho} \mathrm{Var}({Q}_n(b))] \leq 1 + o(1)$ and
\[
\mathrm{Var}(U_n(\rho)) \left\{ \begin{array}{ll} \sim \mathrm{Var}(Q_n(a)) \asymp e^{\rho} \|\eta\|_1^2 \|\mu\|_1^2 \|\nu\|_1^2, &\quad \mbox{if $e^{\rho} \ll e^{\tilde{\rho}}$ (Case (S))},  \\
\asymp e^{\rho} \|\eta\|_1^2 \|\mu\|_1^2 \|\nu\|_1^2, &\quad \mbox{if $e^{\rho} \asymp e^{\tilde{\rho}}$ (Case (L1))}, \\ \sim \mathrm{Var}(e^{\rho} Q_n(b)) \asymp e^{2\rho}(\mu, \eta)(\nu, \eta) \|\mu\|_1^2 \|\nu\|_1^2, &\quad \mbox{if $e^{\rho} \gg e^{\tilde{\rho}}$ (Case (L2))},  \\ \end{array} \right.
\]
where as before $\tilde{\rho} = \tilde{\rho}(\mu, \nu, \eta) = \log(\|\eta\|_1^2/[(\mu, \eta)(\nu, \eta)])$.
\item
We have ${\rm SNR} \asymp 1 / \sqrt{r_n(\theta)} \goto \infty$.
Especially, in the MDH case where (\ref{condition:moderate}) holds, 
\[
\mathrm{SNR}  \asymp   n e^{\gamma/2} \cdot \min\{1, e^{(\rho+ \gamma)/2}\} = \left\{
\begin{array}{ll}
n e^{(\rho/2) + \gamma}, &\qquad \mbox{$e^{\rho} \gg e^{-\gamma}$ \; (Case (S))},  \\
ne^{\gamma/2}, &\qquad e^{\rho} \geq \mbox{$C e^{-\gamma}$  \; (Case (L))}.  \\ 
\end{array}
\right.
\]
\end{itemize}
\end{theorem}
Theorem \ref{thm:snr1} is proved in Section \ref{sec:proof}, with more details in the supplement.
The proof is difficult and extremely long.  One reason is that, $Q_n(a)$, $Q_n(b)$ and
$U_n(\rho)$   are U-statistics in nature. The analysis of such statistics frequently involves
delicate combinatorics, and is known to be hard, tedious, and error-prone. Another reason is that, we do not want to impose any artificial conditions (so we may have a shorter proof, but
 may also hurt the practical utility). See Section \ref{sec:proof} and also \cite{JKL2021, JKLW2022} for more detailed discussion on this.

\subsection{Asymptotical optimality of the LCR estimate $\hat{\rho}^*$}  \label{subsec:est}
Without loss of generality,  we measure the loss by mean-squared error (MSE)  for any estimator for $\rho$. 
Recall that in condition   (\ref{condition2}), we assume $\|\eta\|_1 \goto \infty$
and $e^{\rho/2} \|\eta\|_1 \goto \infty$.  In Theorem \ref{thm:estUB} below, we {\it slightly} strengthen this condition  by assuming
\begin{equation} \label{conditionstr1}
\|\eta\|_1 / \log^2(n) \goto \infty,  \qquad \mbox{and} \qquad  e^{\rho/2} \|\eta\|_1 / \sqrt{\log(n)} \goto \infty. 
\end{equation}
Recall that $\mathrm{SNR} \asymp 1 / \sqrt{r_n(\theta)}$ and note that by definitions and basic algebra,
\begin{equation} \label{Definern}
r_n(\theta)  =  \max\bigl\{1/[e^{\rho} \|\eta\|_1^2], \; (\mu, \eta)(\nu, \eta) /\|\eta\|_1^4\bigr\} \asymp
\left\{
\begin{array}{ll}
1 / [e^{\rho} \|\eta\|_1^2], &\mbox{in Case (S)},  \\
(\mu, \eta) (\nu, \eta) /\|\eta\|_1^4, &\mbox{in Case (L)}.  \\
\end{array}
\right.
\end{equation}
The following theorem is proved in the supplement.
\begin{theorem} \label{thm:estUB}
{\it (Upper bound)}.
Suppose (\ref{condition1}), (\ref{conditionstr1}), and (\ref{condition3})-(\ref{condition4}) hold.  As $n \goto \infty$,
\[
\mathbb{E}[(\hat{\rho}^* - \rho)^2] \leq C \frac{\mathrm{Var}(U_n(\rho))}{(\mathbb{E}[Q_n(a)])^2} \leq   C r_n(\theta),
\]
where $r_n(\theta) \goto 0$.
If additionally (\ref{condition:moderate}) holds,  then $r_n(\theta) \leq (C / n^2 e^{\gamma}) \max\{1, e^{-(\gamma + \rho)}\} \goto 0$.
\end{theorem}
Introduce $c_{\mu, \nu, \eta} = \|\eta\|_1^4 / [(\mu, \eta) (\nu, \eta) \|\mu\|_1 \|\nu\|_1]$ and
\begin{equation} \label{Definern-}
	r_n^{-}(\theta) = \max\bigl\{1/[e^{\rho} \|\eta\|_1^2], \; 1/\|\mu\|_1 \|\nu\|_1\bigr\} \asymp \left\{
	\begin{array}{ll}
		\frac{1}{e^{\rho} \|\eta\|_1^2}, &\quad \mbox{Case (S)}, \\
		\frac{1}{\|\mu\|_1 \|\nu\|_1} = c_{\mu, \nu, \eta} (\mu, \eta)  (\nu, \eta) / \|\eta\|_1^4, &\quad \mbox{Case (L)}.  
	\end{array}
	\right.
\end{equation}
\begin{theorem}  \label{thm:estLB}
	(Lower bound).  Suppose conditions \eqref{condition1}-\eqref{condition2} hold. As $n \goto \infty$,
	\[
	\inf _{\hat{\rho}} \sup _{e^\rho \in\left[n^{-2}, n^2\right]} \mathbb{E}\left[\left(\hat{\rho}-\rho\right)^2\right] \gtrsim r_n^{-}(\theta),
	\]
	where $r_n^-(\theta) \goto 0$. If additionally (\ref{condition:moderate}) holds,  then $r_n^-(\theta) \geq (C / n^2 e^{\gamma}) \max\{c_{\mu,\nu,\eta}, e^{-(\gamma + \rho)}\} \goto 0$.
\end{theorem}
For an intuitive explanation, recall that the effective sample sizes for directed and reciprocal edges are 
$\|\mu\|_1 \|\nu\|_1$ and $e^{\rho}\|\eta\|_1^2$, respectively. Consider the MDH setting for simplicity.  
In this setting, $e^{\rho} \|\eta\|_1^2 \ll \|\mu\|_1 \|\nu\|_1$ in Case (S) where directed edges are more prevalent than reciprocal edges, and $e^{\rho} \|\eta\|_1^2 \gg \|\mu\|_1 \|\nu\|_1$ in Case (L) where reciprocal edges are more prevalent than directed  edges. Intuitively, the lower bound for MSE should be at the order of $1 / \min\{n_1, n_2\}$, 
where $n_1$ and $n_2$ are the effective sample sizes for directed and reciprocal edges. This explains our lower bound in Theorem \ref{thm:estLB}.

Comparing this with $r_n(\theta)$ in (\ref{Definern}), the RHS are the same except for the factor $c_{\mu, \nu, \eta}$.  If
\begin{equation} \label{condition:mu}
c_{\mu, \nu, \eta} \geq c_1,   \qquad \mbox{for a constant $c_1 >  0$},
\end{equation}
then $r_n(\theta) \asymp r_n^{-}(\theta)$ and the lower bound is tight.
Note that condition (\ref{condition:mu}) is much weaker than the MDH condition (\ref{condition:moderate}).  Therefore,  in the minimax framework, if  
(\ref{condition:mu}) holds, then the construction is a least-favorable configuration and our method is asymptotic minimax.

In detail, we consider a minimax setting as follows. Recall that our parameters are $\theta = (\rho, \gamma, \alpha_1, \beta_1, \ldots, \alpha_n, \beta_n)$. Consider the class of parameters
\[
\Theta_0 = \biggl\{\theta: \min\bigl\{\frac{1}{\mu_{max}}, \frac{1}{\nu_{max}}, \frac{1}{e^{\rho/2} \eta_{max}}, \frac{\|\eta\|_1}{\log^2(n)}, \frac{e^{\rho/2} \|\eta\|_1}{\sqrt{\log(n)}},  \frac{(\mu, \eta)(\nu, \eta)}{\mu_{max}^2 \nu_{max}^2} \bigr\} \geq \log(\log(n)) \biggr\},
\]
where the term $\log(\log(n))$ is chosen only for convenience and can be replaced by other diverging sequence $a_n$.
Note that if we neglect the two multi-$\log(n)$ terms on the left hand side, then $\Theta_0$ defines essentially all $\theta$ satisfying conditions (\ref{condition1})-(\ref{condition3}) hold. 
Fix an $\eps_n > 0$ and consider the class of parameters
\[
\Theta(\eps_n)  = \{\theta \in \Theta_0: r_n(\theta) \leq \eps_n \; \mbox{and condition  (\ref{condition4}) holds in Case (L)}\}.
\]  
The minimax risk in $\Theta(\eps_n)$ is then $R_n^*(\eps_n) = \inf_{\hat{\rho}} \sup_{\theta \in \Theta(\eps_n)} \mathbb{E}[(\hat{\rho} - \rho)^2]$, and the worst-case MSE for $\hat{\rho}^*$ is then $R_n(\hat{\rho}^*; \eps_n) =  \sup_{\theta \in \Theta(\eps_n)} \mathbb{E}[(\hat{\rho}^* - \rho)^2]$. By our conditions, the interesting range for $r_n(\theta)$ (and so that for $\eps_n$)  is 
$n^{-2} \ll r_n(\theta) \ll 1$.  
 Theorem \ref{thm:minimax} is proved in the supplement.
 \begin{theorem}  \label{thm:minimax}
 (Optimality).  As $n \goto \infty$, if  $(\log\log(n) /n)^2 \leq \eps_n\leq 1/(\log(n)[\log\log(n)]^2)$, then $C \eps_n \leq R_n^*(\eps_n) \leq R_n(\hat{\rho}^*; \eps_n)  \leq C \eps_n$, so $\hat{\rho}^*$ is asymptotically minimax. 
 \end{theorem}
 In summary, first, the LCR estimator $\hat{\rho}^*$ is asymptotically minimax and is thus rate optimal. second, the testing result in Theorem \ref{thm:estLB} is sharp in Case (S). Recall that
Case (S) includes the case of $\rho \leq C$ and is the most interesting case in
practice.  The testing result is also sharp in Case (L) provided with (\ref{condition:mu}) holds. Note that condition (\ref{condition:mu}) holds in the MDH case where (\ref{condition:moderate}) holds, but it covers a setting much broader than the MDH case.
We will revisit these in Section \ref{subsec:test} when we discuss the phase transition  (similar to the minimax framework, phase transition is a theoretical framework for optimality (e.g., \cite{DJ15}).

\begin{figure}
\centering
\includegraphics[height = 1.5 in]{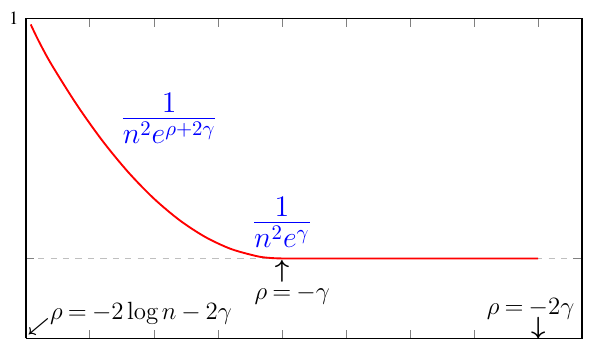}
\caption{Rate of minimax risk ($\gamma$ is fixed);  $x$-axis:  $\rho$ range from $-2(\gamma + \log(n))$ to $-2 \gamma$.}
\label{fig:rate}
\end{figure}

Recall that  the minimax rate is determined by $r_n(\theta)$. It is of interest to consider the MDH case and derive a more explicit form for $r_n(\theta)$.  In the MDH case where (\ref{condition:moderate}) holds, all $\mu_i, \nu_i, \eta_i$ equal  $O(e^{\gamma/2})$ by definitions, and $\tilde{\rho}(\mu, \nu, \eta) = \log(\|\eta\|^2/[(\mu, \eta)(\nu, \eta)]) \asymp - \gamma$. Therefore,
\[
r_n(\theta)  \asymp   \max\{1/(n^2 e^{\rho + 2 \gamma}), 1/(n^2 e^{\gamma})\} \asymp \left\{
\begin{array}{ll}
1/ (n^2 e^{\rho + 2 \gamma)}, &\qquad \mbox{if $e^{\rho} \ll  e^{- \gamma}$ \;  (Case (S))} , \\
1/(n^2 e^{\gamma}), &\qquad \mbox{if $e^{\rho}   \geq C e^{- \gamma} $ \;  (Case (L))}.
\end{array}
\right.
\]
See Figure \ref{fig:rate}. In the MDH case, condition (\ref{condition:mu}) always holds, and our estimator is optimal in the minimax sense (our test also achieves the optimal phase transition; see Section \ref{subsec:test}). 

{\bf Remark 6}.  The $p_1$-model can be viewed as the extension of the $\beta$-model for  directed networks.  The $\beta$-model was well-studied, but existing estimators for the parameters (e.g.,  \cite{Chatterjee2011random, Karwa:Slakovic:2016, stein2020sparse,  Yan:Leng:Zhu:2016, Yan2013clt} are only proved to be consistent when the network density $\tau_n \geq C n^{-1/6}$, where  
$\tau_n =  n^{-2} \times (\mbox{expected number of total edges}$.    For instance, \cite{Chatterjee2011random,Karwa:Slakovic:2016} assumed that  $\tau_n \geq 1/\log n$ or $\tau_n \geq 1/n^{22}$, and \cite{Chen:2020, stein2020sparse} assumed that $\tau_n \gg O(n^{-1/2})$. 
In comparison, we allow all networks with $n^{-2} \log(n) \ll \tau_n \ll 1$, where our estimator $\hat{\rho}^*$ is consistent  $\rho$. 

\subsection{Variance estimation and optimal phase transition for testing}   \label{subsec:test}
We now test $H_0: \rho = \rho_0$, where $\rho_0$ is given but (as before) may depend on $n$. For this purpose,  previous results   are inadequate, for we not only need a good estimator for $\rho$, but also need a consistent estimate for the variance of the estimator. 
Recall that $\hat{\rho}^* = H(\hat{\rho}; 2 \log(n))$ (see (\ref{DefineLCR})) where $\hat{\rho} - \rho \approx U_n(\rho) / \mathbb{E}[Q_n(a)]$
with  $U_n(\rho) = Q_n(a) - e^{\rho} Q_n(b)$ as in (\ref{DefineUn}).  The following lemma is proved in the supplement.
\begin{theorem} \label{thm:normal}
({\it Normality}).
Suppose conditions (\ref{condition1})-(\ref{condition4}) hold.  As $n \goto \infty$,
\[
\frac{U_n(\rho)}{\sqrt{\mathrm{Var}(U_n(\rho))}}  \goto N(0,1)  \;\;
\mbox{and} \;\;    \mathrm{SNR} \cdot (\hat{\rho}^* - \rho) =  \frac{\mathbb{E}[Q_n(a)]}{\sqrt{\mathrm{Var}(U_n(\rho))}}   \cdot (\hat{\rho}^* - \rho) \goto N(0,1),  \;\;
\mbox{in law}.\]
\end{theorem}

We now test the null hypothesis $H_0^{(n)}: \rho = \rho_0$.  By Theorem \ref{thm:normal},  we can tackle this using $(\hat{\rho}^* - \rho_0)$,  but to do so,  we must have a consistent estimate for $\mathrm{SNR}$, or equivalently $\E[Q_n(a)]$ and (especially) $\mathrm{Var}(U_n(\rho))$.    By Theorem \ref{thm:snr1},
\[
\mathrm{Var}(U_n(\rho)) = (1 + o(1))  V_n(\rho), \qquad \mbox{where $V_n(\rho) = \sum_{i \neq j} [2 r_{ij}^2 \Omega_{ij}^{11} + (s_{ij}  - e^{\rho} t_{ij})^2 \Omega_{ij}^{10}]$}.
\]
We now discuss how to estimate $V_n(\rho)$. First, we can estimate $\rho$ by $\hat{\rho}$.  Second,
we can estimate $\Omega_{ij}^{11}$ and $\Omega_{ij}^{10}$ by $A_{ij}^{11}$ and $A_{ij}^{10}$, respectively.
Last, by the sparsity assumption (\ref{condition1}), it is seen that $r_{ij}  = (1 + o(1)) \cdot  \sum_{k, \ell \notin \{i, j\}, k \neq \ell}   \Omega_{k\ell}^{01}$, $s_{ij} = (1 + o(1)) \cdot   \sum_{k, \ell \notin \{i, j\}, k \neq \ell}  \Omega_{k\ell}^{11}$, and $t_{ij}  = (1 + o(1)) \cdot  \sum_{k, \ell \notin \{i, j\}, k \neq \ell} (\Omega_{jk}^{10} \Omega_{k\ell}^{10}   + \Omega_{k\ell}^{10} \Omega_{\ell i}^{10} + \Omega_{jk}^{10}  \Omega_{\ell i}^{10})$, where $o(1) \goto 0$ uniformly for all $i, j$. Therefore, we can estimate
$r_{ij}$, $s_{ij}$, and $t_{ij}$ by $\hat{r}_{ij}$, $\hat{s}_{ij}$, and $\hat{t}_{ij}$, respectively, where
 \[
\hat{r}_{ij} = \sum_{k, \ell \notin \{i, j\}, k \neq \ell}  A_{k\ell}^{01}, \;\;\;    \hat{s}_{ij}  =   \sum_{k, \ell \notin \{i, j\}, k \neq \ell}  A_{k\ell}^{11},  \;\;\;  \hat{t}_{ij} = \sum_{k, \ell \notin \{i, j\}, k \neq \ell} (A_{jk}^{10} A_{k\ell}^{10}   + A_{k\ell}^{10} A_{\ell i}^{10} + A_{jk}^{10}  A_{\ell i}^{10}).
\]
A convenient estimate for $V_n(\rho)$ is therefore
\begin{equation} \label{DefinehatV}
\widehat{V}_n(\hat{\rho})  = \sum_{i \neq j} [2 \hat{r}_{ij}^2 A_{ij}^{11} + (\hat{s}_{ij}  - e^{\hat{\rho}} \hat{t}_{ij})^2 A_{ij}^{10}]. 
\end{equation}
Combining this with Theorem \ref{thm:normal}, we propose the following test statistics:
 \begin{equation} \label{Definepsin}
 \psi_n^* = \widehat{\mathrm{SNR}} \cdot (\hat{\rho}^* - \rho_0),  \qquad  \mbox{where $\widehat{\mathrm{SNR}} = Q_n(a) / \sqrt{\widehat{V}_n(\hat{\rho})}$}.
\end{equation}
Now, if $\widehat{\mathrm{SNR}}  \approx \mathrm{SNR}$ as expected,    then by elementary algebra and Theorem \ref{thm:normal},
$\psi_n^* \approx  \mathrm{SNR} \cdot (\rho - \rho_0)  + \mathrm{SNR} \cdot (\hat{\rho}^* - \rho)$, where $\mathrm{SNR} 
\cdot (\hat{\rho}^* - \rho) \goto N(0,1)$ by Theorem \ref{thm:normal}.  
Therefore,  the power of the test hinges on the quantity of $\mathrm{SNR} \cdot |\rho - \rho_0|$. 
In detail, fix $0 < \alpha < 1$ and let $z_{\alpha/2}$ be the number such that $\mathbb{P}(N(0,1) \geq z_{\alpha/2}) = \alpha/2$.  Consider the  LCR  test where we reject the null  if and only if
$|\psi_n^*| \geq z_{\alpha/2}$.
 We assume
\begin{equation} \label{conditionstr2}
e^{\rho/2} \|\eta\|_1 / \sqrt{\log(n)} \goto \infty,  \;\;\; \mbox{and} \;\;\; \mbox{$\min\{\frac{\mu_{min}}{\mu_{max}},
\frac{\nu_{min}}{\nu_{max}}\} \cdot \|\eta\|_1 / \log(n) \goto \infty$ in Case (L)}.
\end{equation}
Here, the first condition is mild and appeared before in (\ref{conditionstr1}). The second condition is required in deriving  the  large-deviation bounds for $\widehat{V}_n(\rho)$. The condition is also mild (e.g., in the MDH case, it reduces to
$\|\eta\|_1 / \log(n)  \goto \infty$).  Theorem \ref{thm:test}  is proved in the supplement.
\begin{theorem} \label{thm:test}
Suppose conditions (\ref{condition1}),  (\ref{condition4}) 
and (\ref{conditionstr2}) hold, and consider the two-sided LCR test above for testing the null hypothesis $H_0^{(n)}: \rho = \rho_0$.
As $n \goto \infty$,
\begin{itemize}
\item $\widehat{V}_n(\hat{\rho})/ V_n(\rho) \goto 1$  and $\frac{\widehat{\mathrm{SNR}}}{\mathrm{SNR}}  \goto 1$ in probability, and $\widehat{\mathrm{SNR}}  \cdot  (\hat{\rho}^* - \rho)  \goto N(0,1)$ in law.
\item If $H_0^{(n)}$ holds, then $\psi_n^*  \goto N(0,1)$ in law and the level of the LCR test is $\alpha + o(1)$.
\item If $H_0^{(n)}$ is not true and $\mathrm{SNR} \cdot |\rho - \rho_0|  \goto \infty$, then the power of the LCR test $\goto 1$.
\end{itemize}
\end{theorem}
In the case where $\mathrm{SNR}  \cdot (\rho- \rho_0) \goto c$ for a constant $c$, we can similarly show that the power is $1-\Phi(z_{\alpha/2}-c)+\Phi(-z_{\alpha/2}-c) + o(1)$.  Note that (\ref{condition4}) and (\ref{conditionstr2}) are only required in Case (L)  where $e^{\rho}$ is relatively large.  In Case (S),  we have a much simpler approximation for $\mathrm{Var}(U_n(\rho))$:    $\mathrm{Var}(U_n(\rho)) \sim \sum_{i \neq j} 2 r_{ij}^2 \Omega_{ij}^{11}$, and so we have a simpler estimator for $\mathrm{Var}(U_n(\rho))$ (and so a simpler test).    In detail, denote
$W_n$ by $W_n = \sum_{i \neq j} 2 r_{ij}^2 \Omega_{ij}^{11}$ and let $\widehat{W}_n = \sum_{i \neq j} 2 \hat{r}_{ij}^2 A_{ij}^{11}$. Introduce a simpler version of $\psi_n^*$ by $\phi_n^* = \widehat{\mathrm{snr}} \cdot (\hat{\rho} - \rho_0)$ where $\widehat{\mathrm{snr}}= Q_n(a) / \sqrt{\widehat{W}_n}$. Similarly,
fixing a $\alpha \in (0,1)$, we can have a two-sided test that rejects the null when $|\phi_n^*| \geq z_{\alpha/2}$. We call this the {\it simplified LCR  test}. 
Corollary \ref{cor:test} is proved in the supplement.
\begin{cor} \label{cor:test}
Suppose (\ref{condition1}) and (\ref{conditionstr1}) hold. Consider the simplified  LCR test for Case (S) where $e^{\rho} \cdot (\mu, \eta)(\nu, \eta) /  \|\eta\|_1^2 \goto 0$. As $n \goto \infty$,
\begin{itemize}
\item $\widehat{W}_n / W_n \goto 1$ and $\frac{\widehat{\mathrm{snr}}}{\mathrm{SNR}} \goto 1$ in probability. Also,  $\widehat{\mathrm{snr}} \cdot (\hat{\rho}^* - \rho)  \goto N(0,1)$ in law.
\item If $H_0^{(n)}$ holds, then $\phi_n^*  \goto N(0,1)$ and the level of the simplified LCR test is $\alpha + o(1)$.
\item If $H_0^{(n)}$ is not true and $\mathrm{SNR} \cdot |\rho - \rho_0|  \goto \infty$, then the power of the test above $\goto 1$.
\end{itemize}
\end{cor}
We can slightly relax (\ref{conditionstr1}) by replacing $\|\eta\|_1 / \log^2(n) \goto \infty$ with that of $(\|\mu\|_1\|\nu\|_1)/\log(n)\goto \infty$.    Also, recall that (\ref{conditionstr1}) is a slightly stronger version of (\ref{condition2}),  
so essentially we only require (\ref{condition1})-(\ref{condition2}) here, 
which are very relaxed conditions.    As discussed before, Case (S) covers the most interesting case from a practical viewpoint (especially, it includes the case of $\rho \leq C$). When we believe $\rho$ is relatively small,
we prefer $\phi_n^*$ to $\psi_n^*$, for $\phi_n^*$ is not only simper, but the regularity conditions
in Corollary \ref{cor:test} are more relaxed than those in Theorem \ref{thm:test}.

{\bf Remark 7}. In Theorem \ref{thm:test} and Corollary \ref{cor:test}, the targeted Type I error is $\alpha \in (0,1)$ (fixed).  We may allow $\alpha$ to tend $0$, so the Type I error $\goto 0$:
we call $\alpha = \alpha_n$ tends to $0$ {\it slowly enough}  if $\alpha \goto 0$ and $z_{\alpha/2} \leq C \log(\log(n))$.
Here, the $\log\log(n)$ can be replaced by any sequence that tends  to $\infty$ slowly enough.
In such a case, in Theorem \ref{thm:test} and Corollary \ref{cor:test}, the Type II error of the tests $\goto 0$ provided with a slightly strong condition: $[\log\log(n)]^{-1}  \cdot |\rho - \rho_0| \cdot \mathrm{SNR} \goto \infty$.

We now discuss the lower bound.  Recall that $r_n^{-}(\theta)$ is defined in (\ref{Definern-}).
The following corollary is a direct result of Theorem \ref{thm:estLB} so the proof is omitted.
\begin{cor} \label{cor:LB}
({\it Lower bound}).  Suppose conditions (\ref{condition1})-(\ref{condition2})  hold.
If $|\rho - \rho_0| / \sqrt{r_n^{-}(\theta)}  \goto 0$, then we can pair the alternative with a null, so that
the $\chi^2$-distance between two distributions (null and alternative) $\goto 0$. As a result, the sum of Type I and Type II error of any test $\goto 1$ (so any test will be powerless in separating the alternative from the null).
\end{cor}

Our results describe a sharp phase transition.  Similar as 
minimax,  phase transition is a theoretical framework for optimality. 
However, in the minimax framework, we are primarily focused on
optimality near least favorable cases, but least favorable cases are not
necessarily the most interesting or challenging cases in practice.
The phase transition framework allows us
to consider optimality at {\it all cases}, so practically it can be a more appealing (see \cite{DJ15, phase}).

Recall that in Section \ref{subsec:est}, we show the LCR estimator
$\hat{\rho}^*$ is optimal in minimax sense. Below, we further show
the LCR tests $\psi_n^*$  and $\phi_n^*$ achieve the optimal phase transition. 
In detail, recall that $r_n(\theta) \asymp r_n^{-}(\theta)$ in Case (S) and
$r_n(\theta) \asymp r_n^{-}(\theta)$ in Case (L) provided with (\ref{condition:mu}) holds.
Recall $\tilde{\rho} = \tilde{\rho}(\mu, \nu, \eta) = \log(\|\eta\|_1^2 / [(\mu, \eta)(\nu, \eta)])$.
Comparing Theorem \ref{thm:test}, Corollaries \ref{cor:test}-\ref{cor:LB}, and Remark 7,
we have two case, Case (S) and Case (L).  In Case (S), $e^{\rho} \ll e^{\tilde{\rho}}$ (where $\rho$ is small and moderately large; this is the most interesting case in practice). In this case,   $r_n(\theta) \asymp r_n^{-}(\theta)$. If $|\rho - \rho_0| \ge \log(\log(n)) \cdot \sqrt{r_n(\theta)}$,
then the sum of Type I and Type II error of the LCR tests $\psi_n^*$ and $\phi_n^*$ tend to $0$.
If $|\rho - \rho_0| \ll \sqrt{r_n(\theta)}$, then for all tests,  the sum of Type I and Type II error tends to $1$ (and so all tests fail).
In Case (L). $e^{\rho} \ge C e^{\tilde{\rho}}$ (where $\rho$ is relatively large).
In this case, if $|\rho - \rho_0| \ge \log(\log(n)) \cdot \sqrt{r_n(\theta)}$,
then the sum of Type I and Type II error of the LCR test $\psi_n^*$ tends to $0$.
If $|\rho - \rho_0| \ll \sqrt{r_n(\theta)}$ and condition (\ref{condition:mu}) holds, then $r_n(\theta) \asymp r_n^{-}(\theta)$, and  for all tests,  the sum of Type I and Type II error tends to $1$. 
Therefore, the LCR tests achieve the optimal phase transition in Case (S), and also in Case (L) if (\ref{condition:mu}) holds.  Especially, in the MDH case where (\ref{condition:moderate}) holds, condition  (\ref{condition:mu}) holds automatically. In this case, the LCR test $\psi_n^*$ achieves the optimal transition in both Case (S) and Case (L), and it is optimal in the whole range of interest. Note that the optimality here is for each individual parameter $\theta$, but the optimality in the minimax framework
is essentially only for least favorable cases. In this sense, LCR meets a higher standard of optimality.

\subsection{Comparison with the other two pairs of $(a, b)$ in Figure \ref{tab:pairs}}
\label{subsec:otherpairs}
Recall that the key of LCR is to find a pair of $(a, b) \in S_m$ (i.e., generalized cycles) to satisfy (\ref{critical}).
When $m = 4$, we have identified  three non-isomorphic pairs of $(a, b)$ satisfying  (\ref{critical}), denoted by $(a^{(k)}, b^{(k)})$ for
$k = 1, 2, 3$; see Figure  \ref{tab:pairs}. So far,  we have only considered the first pair.  For the other  two  pairs,  we can similarly construct an estimator  $\hat{\rho}^*$ and a test of $\psi_n^*$. By previous sections,  to compare the performances of $(\hat{\rho}_n^*, \psi_n^*)$
constructed for different pairs of $(a, b)$, it is sufficient to compare the SNR.
In light of this, similarly as in (\ref{DefineSNR}), for $k = 1, 2, 3$, define the SNR for the resultant estimator and test associated with $(a^{(k)}, b^{(k)})$ by
\[
\mathrm{SNR}_k = \mathbb{E}[Q_n(a^{(k)})] / \sqrt{\mathrm{Var}(U_{n, k}(\rho))}, \qquad
\mbox{where\; $U_{n, k}(\rho) =  Q_n(a^{(k)}) - e^{c_k  \rho} Q_n(b^{(k)})$},
\]
and for short,  $c_k = c_0(a^{(k)}, b^{(k)}) = 1, 2, 1$ for $k = 1, 2, 3$. Let 
$c_{\mu} = [\|\eta\|_1 (\mu \circ \mu, \eta)]/ [(\mu, \eta)^2]$ and $c_{\nu}
= [\|\eta\|_1 (\nu \circ \nu, \eta)] / [(\nu, \eta)^2]$.  
\begin{table}[htb!]
\renewcommand\arraystretch{1}
\scalebox{0.82}{
\begin{tabular}{lll}
\hline
$k$   & The general case & The MDH case   \\
\hline
1   & $\mbox{SNR}_1 \asymp \|\eta\|_1^2 / \sqrt{\max\{e^{-\rho} \|\eta\|_1^2, (\mu, \eta) (\nu, \eta)\}}$ & $n e^{\gamma/2} \min\{1, e^{(\rho + \gamma)/2}\}$   \\
2  &  $\mbox{SNR}_2 \asymp    \|\eta\|_1^2 / \sqrt{\max\{e^{-\rho} \|\eta\|_1^2, (\mu, \eta) (\nu, \eta), 1\}}$ & $n e^{\gamma/2} \min\{1, e^{(\rho + \gamma)/2}, e^{3 \gamma/2}\}$ \\
3 & $\mbox{SNR}_3 \leq C \|\eta\|_1^2 / \sqrt{\max\{e^{-\rho} \|\eta\|_1^2 \max\{c_{\mu}, c_{\nu},1/[(\mu, \eta) (\nu, \eta)]\}, (\mu, \eta) (\nu, \eta)\}}$ & $n e^{\gamma/2} \min\{1, e^{(\rho + \gamma)/2}, ne^{(\rho+4\gamma)/2}\}$ \\
\hline
\end{tabular}
}
\caption{SNR comparison, each is for an $\hat{\rho}^*$ corresponding to a pair of  $(a, b)$  in  Figure \ref{tab:pairs}.}
\label{tab:rates}
\end{table}
Here,  (a) $c_{\mu}, c_{\nu} \geq 1$ by Cauchy-Schwartz inequality, and (b) while
we may have $\max\{c_{\mu}, c_{\nu}\} \gg 1$ in the SDH case,  $\max\{c_{\mu}, c_{\nu}\} \leq C$   in the  MDH  case where (\ref{condition:moderate}) holds.
Table \ref{tab:rates} compares  the three SNR, with more details in Table A4-A5 of the supplement.

In summary, we have the following. 
First, among three SNR's, $\mathrm{SNR}_1$ is always the largest and the corresponding estimator is 
optimal. Second, the order of $\mathrm{SNR}_3$  may be
significantly smaller than that of $\mathrm{SNR}_1$, but is the same if $\max\{c_{\mu}, c_{\nu},1/[(\mu, \eta) (\nu, \eta)]\} \leq C$.  
Therefore, the corresponding estimator is optimal in the subregion of the parameters where $\max\{c_{\mu}, c_{\nu},1/[(\mu, \eta) (\nu, \eta)]\} \leq C$,  but may be non-optimal outside the region.    
Similarly, the order of $\mathrm{SNR}_2$  is the same as that of $\mathrm{SNR}_1$  if $\max\{e^{-\rho} \|\eta\|_1^2, (\mu, \eta) (\nu, \eta)\}  \geq C$ and is significantly smaller if $\max\{e^{-\rho} \|\eta\|_1^2, (\mu, \eta) (\nu, \eta)\}  \ll 1$. 
Therefore, the corresponding estimator is optimal in the subregion where $\max\{e^{-\rho} \|\eta\|_1^2, (\mu, \eta) (\nu, \eta)\}  \geq C$ but may be non-optimal outside the region.    Take the
MDH case for example.
We can divide the parameters $(\rho, \gamma)$ into two cases:
Case 1:   $\gamma  > -(2/3) \log(n)$, and
Case 2:   $\gamma  \le -(2/3) \log(n)$. 
Both the orders of $\mathrm{SNR}_2$ and $\mathrm{SNR}_3$ are the same as that of $\mathrm{SNR}_1$ in Case 1,
but can be significantly smaller in Case 2 when $\rho$ is of the order $O(1)$. 
Therefore,  all three estimators are optimal in Case 1, but only the first one 
is optimal in both cases.   For these reasons, we prefer to use the first pair of $(a, b)$ in Figure \ref{tab:pairs} to construct the LCR estimators and tests. Note that for each of three pairs of cycles in Figure 1, we construct an estimator for $\rho$,  
and in the noiseless case, all estimators equal to $\rho$.    Consider now the real case. 
Note that the expected numbers of cycles in the first pair are larger than those 
in the other two pairs. Therefore, the first estimator has more effective samples 
to estimate $\rho$ and so has smaller variance. 
 This explains why we have the best rate if we choose the first pair in Figure \ref{tab:pairs}  to construct 
our statistic.

\subsection{Extensions and future work}  \label{subsec:extension}
The $p_1$ model has $(2n+2)$ parameters $\rho, \gamma, \alpha_1, \beta_1, \ldots, \alpha_n, \beta_n$. For
reasons of space and especially practical interests, we focus on the reciprocal parameter $\rho$ in this paper, but our idea is readily extendable to estimating $\gamma$ (overall sparsity) and $(\alpha_i,\beta_i)$ (heterogeneity parameter for node $i$).
In fact, by carefully picking an $m \geq 3$ and two types of generalized $m$-cycles,  Type $C$ and Type $D$, we discover that for all $1 \leq i \leq n$,  we have $\rho + \gamma + \alpha_i  + \beta_i = \log(C_i / D_i)$,  where $C_i$ and $D_i$ are the expected
number cycles containing node $i$ of Type C and Type D, respectively.  Such a discovery can be used to
develop a convenient estimate for all other $(2n+1)$ parameters of $p_1$ model.  Here, similarly as in Section \ref{subsec:trick},
each $C_i$ and $D_i$ is a big sum, containing the nonlinear factors $K_{ij}$ of the $p_1$ model, but if we design the
two types of $m$-cycles carefully, the nonlinear term cancel out with each other in the ratio of $C_i / D_i$.

Note the idea we just mentioned is readily extendable to the $\beta$-model. In the $\beta$-model, we similarly have
$A = \Omega - \diag(\Omega) + W$ as in (\ref{plowrankmodel1}) but where $(A, \Omega, W)$ are symmetrical
and $\log(\Omega_{ij} / (1 - \Omega_{ij})) = \beta_i + \beta_j$.  How to estimate $\beta_1, \ldots, \beta_n$ in
the $\beta$-model is also a long-standing problem (e.g., \cite{Chatterjee2011random, Chen:2020, Yan:Leng:Zhu:2016, Yan2013clt}),
but fortunately, since the $\beta$-model can be viewed as a special case of the $p_1$ model,
the idea above is extendable to address the problem.
Of course, while the high-level idea is extendable, we still need delicate (and especially very long)
analysis to study the estimators for the new settings. Since this paper is already very long (the supplement is about $150$  pages), we leave the study along this line to the future.

\begin{table}[!htp]
\centering
\caption{Comparison with $|\hat{\rho}^*-\rho|$ (left) and $|\hat{\rho}^{mle} - \rho|$ (right); based on $100$ repetitions.
NA:  MLE failed to exist in all $100$ repetitions}
\label{table-compare}
\scalebox{0.8}{
\begin{tabular}{ccccccc}
$\gamma$ & $n$ & $\rho=-\tfrac{1}{4}\log n$  & $\rho=0$ & $\rho=0.5$ & $\rho=\log(\log n)$   &  $\rho=\tfrac{1}{4}\log n $    \\ 
\hline 
$\gamma_1$  & $500$  & $(.024, .025)$&$ (.018, .018)  $&$ (.016, .016)   $&$ (.016, .016)  $&$ (.015, .015) $ 
\\
            & $1000$ & $ (.015, .014)  $&$  (.008, .008)  $&$ (.008, .008)  $&$ (.007, .008)   $&$ (.007, .008)  $
\\
            & $2000$ & $ (.009, .009)  $&$ (.005, .005)  $&$ (.004, .004)  $&$ (.003, .004)  $&$ (.004, .004)$
                       \\
$\gamma_2$  & $500$    & $ (.071, .069) $&$ (.038, .036) $&$ (.033, .033) $&$ (.023, .022)  $&$ (.023, .024)  $
\\
            & $1000$   & $ (.050, .057) $&$ (.023, 0.021)  $&$ (.020, .019) $&$ (.014, .014)  $&$ (.014, .014)  $
\\
            & $2000$   & $ (.037, .038)  $&$ (.016, .016)  $&$ (.014, .014) $&$ (.008, .008)  $&$ (.009, .009)$
\\
$\gamma_3$  & $500$    & $ (.234,  NA) $&$  (.126, NA) $&$ (.099,  NA)  $&$ (.062, NA) $&$ (.058, NA) $
\\
            & $1000$   & $ (.263, NA)  $&$ (.089, NA) $&$ (.082, NA)  $&$ (.046, NA)  $&$ (.048, NA)  $
\\
            & $2000$   & $ (.201, NA)  $&$ (.091, NA)  $&$ (.061, NA) $&$ (.035, NA)  $&$ (.037, NA) $
\\
\end{tabular}
}
\label{table:compare}
\end{table}


\section{Numerical study} \label{sec:numeric}
\subsection{Simulations}
Our simulated work contains two parts, an estimation part and a testing part.

{\bf Estimation part}.  We generate the adjacency matrix $A$ in the $p_1$ model as follows. For each $n = 500, 1000, 2000$,    we consider  $3$ different choices of $\gamma$ ($\gamma_1 = -\log(n)/4, \gamma_2 =  - \log(n)/2,  \gamma_3 = -\log n + \log(\log n)$)  and $5$ different choices of $\rho$
($\rho_1 = -(\log n)/4$, $\rho_2 = 0$, $\rho_3 = 0.5$, $\rho_4 = \log(\log n)$ and $\rho_5 = (\log n)/4$).
Generate $\alpha_i$ and $\beta_i$ independently from $N(0,1)$ and $U(-1, 1)$, respectively (centered empirically  so that (\ref{p1model2}) holds).
These cover a wide range of  reciprocal effects and network sparsity where the network densities (defined as $n^{-2} \times \mbox{expected number of total edges}$) are   $n^{-1/4}$, $n^{-1/2}$ and $\log n/n$ for $\gamma = \gamma_1, \gamma_2, \gamma_3$, respectively,  and allow relatively severe degree heterogeneity where
the largest degree can be hundred times larger than the smallest one. 

\begin{figure}[!htp]
\centering
\includegraphics[height=1.5 in, width=5.05in, angle=0]{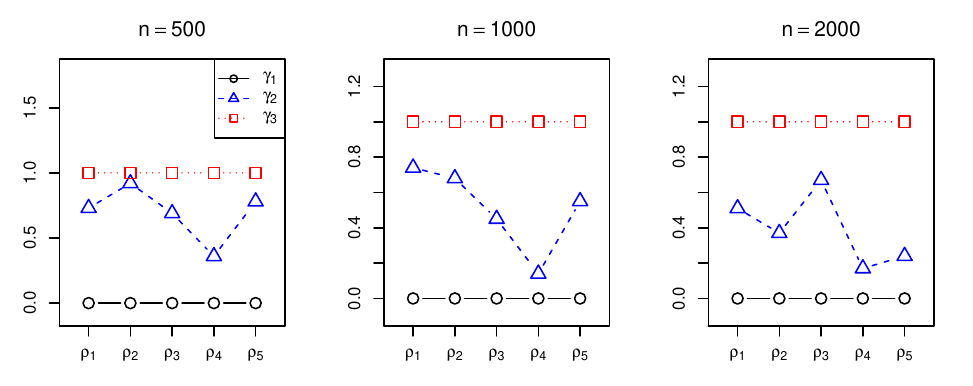} 
\caption{The fraction (out of $100$ repetitions) that the MLE does not  exist in each setting.}
\label{mle-fail}
\end{figure}

First, we compare our estimator $\hat{\rho}^*$ with the MLE estimator $\hat{\rho}^{mle}$.  
We implemented the MLE using the frequency iterative algorithm (FIA) proposed by \cite{Holland:Leinhardt:1981} in R on a computer with an Intel T7700 2.40GHz processor and a 16GB RAM.
Note that in order for MLE to exist, we need  strong (and generally hard-to-check) conditions
\cite{Rinaldo2013mle}. For example, when there is a node whose in-degree or out-degree is $0$ or $(n-1)$, the MLE does not exist (the maximum likelihood is $\infty$). For each setting, we repeat $100$ times, but for MLE, we only count the repeats where the MLE exists (this is in favor of MLE).  The results are in Table \ref{table-compare} (NA:  MLE does not exist for all $100$
repetitions).  For each setting, the fraction that the MLE does not exist (in $100$ repetitions) is reported in
Figure \ref{mle-fail}. Computation time (in R) of two approaches are  in  Table \ref{table-compare-time}.
 We observe that, first,  our estimator LCR generally perform satisfactorily, especially for larger
$(n, \rho)$ and smaller $\gamma$ (where it is more sparse). Second,  the MLE fails to exist in many cases, but
when it exists, it also works well, with an accuracy comparable to the LCR.
Last, in the experiments we investigated, the MLE is slower than our methods by $6$-$16$ times in computation.
Therefore, the LCR has a clear winning edge over the MLE.

Our approach is relatively robust to model misspecification (e.g., 
contaminated nodes, misspecified number 
of communities). Consider a misspecified setting where we have a directed network with $2$ communities.  For a  parameter $\theta$,  assume $\mathbb{P}(A_{ij} = a, A_{ji} = b) =
K_{ij} \cdot \mathrm{exp}(a (\gamma + \alpha_i + \beta_j) + b (\gamma + \alpha_j + \beta_i) + a b \rho  + (a+b) \theta \delta_{ij})$, where $
K_{ij} = [1 + e^{\gamma + \alpha_i + \beta_j + \theta \delta_{ij}}  + e^{\gamma + \alpha_j + \beta_i + \theta \delta_{ij}}  +
e^{2\gamma + \alpha_i + \beta_j + \alpha_j + \beta_i + \rho+ 2\theta  \delta_{ij}}]^{-1}$ and $\delta_{ij} = 1$ if nodes $i, j$ belong to the same community and $\delta_{ij} = 0$ otherwise.   
In such a setting, LCR may be biased, but for a wide range of $\theta$, the biases are 
relatively small. See Table \ref{Table:b} below, where $n = 1000$ and $\gamma =  -\log(n)/4$.  
This suggests that our procedure is relatively robust to model misspecification.


	\begin{table}[htb!]
	\centering
	\caption{Computing time for our method LCR (left) and MLE (right) in seconds (average in 10 repetitions). In many cases, MLE is $10$ times slower than our method ($ \gamma_1 = - \log(n)/4$).}
	\label{table-compare-time}
	\scalebox{0.85}{
		\begin{tabular}{ccccccc}
			\hline
			$\gamma$ & $n$   & $\rho=-\tfrac{1}{4}\log n$  & $\rho=0$ & $\rho=0.5$ & $\rho=\log(\log n)$   &  $\rho=\tfrac{1}{4}\log n $      \\
			\hline
			\multirow{3}{*}{$\gamma_1$} & $1000$  & $1.843\quad 26.341  $&$1.895\quad 22.075  $&$ 1.883\quad 25.493 $&$ 1.994\quad 30.974 $&$ 1.943\quad 24.749 $ \\
			& $2000$ & $9.28\quad 80.156 $&$ 9.766\quad 98.839 $&$ 9.826\quad 100.785  $&$ 10.787\quad 102.422 $&$ 10.691\quad 146.662  $ \\
			& $3000$ & $27.486\quad 170.712  $& $29.169\quad 290.237 $& $ 29.512\quad 188.213 $& $ 32.192\quad 249.429  $& $ 31.684\quad 232.930 $ \\
			\multirow{3}{*}{$0.5$} & $1000$  & $1.908\quad 41.116 $&$1.538\quad 39.549  $&$ 1.579\quad 41.158 $&$ 1.508\quad 51.229 $&$ 1.519\quad 50.109 $ \\
			& $2000$ & $9.824\quad 201.371 $&$ 8.302\quad 186.122 $&$ 8.66\quad 184.916  $&$ 8.092\quad 186.748 $&$ 8.143\quad 208.903  $ \\
			& $3000$ & $29.978\quad > 300  $& $28.565\quad > 300 $& $ 31.18\quad > 300 $& $ 30.989\quad > 300  $& $ 31.081\quad > 300 $  \\
			\hline
		\end{tabular}
	}
\end{table}
	
\begin{table}[!htp]\centering
\caption{Average absolute biases when the model is misspecified ($\theta\neq 0$).}
\label{Table:b}
\vskip5pt
\begin{tabular}{cccc cc}
\hline
$\rho$   & $\theta=-0.5$ & $\theta=-0.2$  & $\theta=0$  & $\theta=0.2$ & $\theta=0.5$   \\
\hline
$0.5$   & 0.062 & 0.012 & 0.008 & 0.012 & 0.061     \\
        $1$     & 0.061 & 0.013 & 0.007 & 0.011 & 0.061  \\
\hline
\end{tabular}
\end{table}

\begin{figure}[!htp]
\centering
\includegraphics[height=1.5 in, width=5.25in, angle=0]{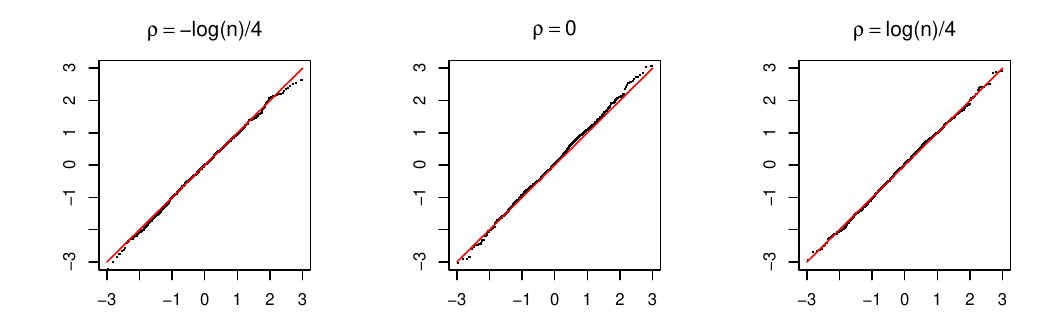}
\caption{Q-Q plots for $\psi_n^*$ when the null is true ($n = 1000$,  $\gamma = - \log(n)/2$).}
\label{fig-rho}
\end{figure}

{\bf Testing part}. We now discuss our test $\psi_n^*$.  First,  Theorem \ref{thm:test} shows
that the limiting null of $\psi_n^*$  is $N(0,1)$.   Figure \ref{fig-rho}  presents the Q-Q plots of
simulated $\psi_n^*$ when the null hypothesis is true ($n = 1000$,  $\gamma = \gamma_2 = - \log(n)/2$, $\rho = \rho_1, \rho_2, \rho_5$, repeat for $1000$ times). The result fits well with our theory, even when the network is relatively sparse.


\begin{figure}[!htpb]
	\centering
	\includegraphics[height = 2.25 in, width=6in]{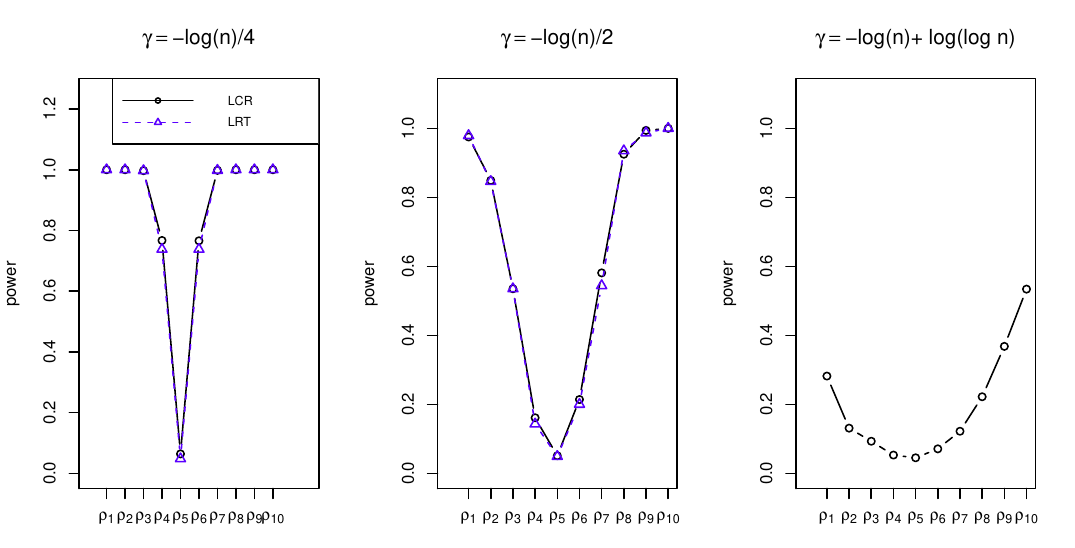}
	\caption{Powers of LCR and LRT ($n=500$, $1000$ repetitions).}
	\label{fig:power}
\end{figure}

Second, we study the power.  There are relatively few testing
approaches for the $p_1$ model, among which there is the
(generalized) LRT suggested by Holland and Leinhardt \cite{Holland}. Since LRT
is closely related to MLE, it is not surprising that (a) algorithm for LRT is only available for $\rho_0 = 0$ (so we take $\rho_0$ below), (b) LRT does not exists in many cases,
and (c) even when it exists, the underlying theory is
largely unknown, but it was conjectured to have limiting null of $\chi_1^2(0)$ \cite{Holland}.
We now compare the powers  of two tests at the nominal level of $5\%$ (using $N(0,1)$ and $\chi_1^2(0)$ as
the null distribution, respectively).  For the alternative, take $\rho =-0.2,-0.15, -0.1, \dots,0.25$ but other parameters are the same as the null. We repeat $1000$ times (and to save time, we only report that for $n = 500$). The results are in Table \ref{table-coverage} and Figure \ref{fig:power}.
The powers of the LCR method and the LRT are very close when $\gamma=-\log n/4, -\log n/2$ while the MLE failed to exist in all simulations when $\gamma=-\log n +\log\log n$,
corresponding to a extremely sparse network with network densities $0.02$.
Note that when $\rho = 0$ (first column), we are in the null case, so the power is nothing but the Type I error; since the target level is $5\%$, the results fit well with our theory. Also, when $\rho$ gradually deviate from $\rho_0$, the power
rapidly increase to $1$, as expected (the last row of Table \ref{table-coverage}  corresponds to the case of $\gamma_3 = -\log(n) + \log\log(n)$,    where the networks are very sparse (network density is $\log(n)/n$). In this case,  for all $\rho$'s  in Table \ref{table-coverage}, the SNR are relatively low, so the powers are also relatively low. We can however increase the SNR and power by 
increasing $\rho$. For example, when $\rho = 0.5$,  the power is $0.98$). Note that LRT does not exist for many cases, especially for sparse networks (i.e., a smaller $\gamma$). Note also LRT is computationally much slower than our test.
These results suggest that our test has a clear winning edge of the LRT. 


\begin{table}[htb!]
	\centering
	\caption{Powers of our methods (first cell value) and LRT (second cell value); cell values are power $\times 100$ ($n = 500$, $1000$ repetitions). Boldface: fraction of repetitions where LRT does not exist (NA: LRT does not exist for all  repetitions; note that $\rho = 0$ is the null case as $\rho_0 = 0$.}
	\label{table-coverage}
	\scalebox{0.5}{
		\begin{tabular}{ccccccccccc}
			& $\rho=-0.20$ & $\rho=-0.15$ & $\rho=-0.10$  & $-0.05$ & $\rho=0$  & $\rho=0.05$ & $\rho=0.10$ & $\rho=0.15$ & $\rho=0.20$ & $\rho=0.25$  \\
			$\gamma_1$  &      $(100, 100,  {\bf 0})$ & $(100, 100, {\bf 0}) $&$(99.7, 99.7, {\bf 0})$ &$(76.7, 73.8, {\bf 0.1}) $ &      $(6.4, 4.8,  {\bf 0.7})$ & $(76.6, 73.8, {\bf 0}) $&$(99.8, 99.8, {\bf 0.1})$ &$(100, 100, {\bf 0.6}) $ &      $(100, 100,  {\bf 0})$ & $(100, 100, {\bf 0.6}) $ \\
			$\gamma_2$  &    $(97.5, 97.9, {\bf 81.3})$ & $(84.9, 84.6, {\bf 70.8 }) $ & $(53.5, 53.5,  {\bf  65.4})$ &
			$(16.1, 14.3, {\bf 79})$ & $(5.1, 4.9, {\bf 81.8})$ & $(21.4, 20, {\bf 56.9})$ & $(58.1, 54.4, {\bf 74.8})$ & $(92.5, 93.5, {\bf 70.6})$ & $(99.4, 98.7, {\bf 69.7})$ & $(100, 100, \bf{75.5})$  \\
			$\gamma_3$  & $(28.2, \mbox{NA}, {\bf 100}$) &  $(13.1, \mbox{NA}, {\bf 100})$ & ($9.3, \mbox{NA}, {\bf 100})$ & $(5.3, \mbox{NA}, {\bf 100})$ & ($4.5, \mbox{NA}, {\bf 100})$ & $(7.1, \mbox{NA}, {\bf 100})$ & ($12.2, \mbox{NA}, {\bf 100})$ & $(22.2, \mbox{NA}, {\bf 100})$ & ($36.8, \mbox{NA}, {\bf 100})$ & $(53.4, \mbox{NA}, {\bf 100})$ \\
		\end{tabular}
	}
\end{table}

\subsection{Analysis of several real networks}
\label{subsec:real} 
We apply the proposed estimator and tests to the well-known weblog data and
a citation network of statisticians.  These networks have more than one communities. 
Since the $p_1$ model is for networks  with only one community,  
it is not appropriate to use it to directly model the entire network, but we can divide each network into several sub-networks, each with only one community.  Now, for each subnetwork, it is appropriate to model it with the $p_1$ model, where our  approaches are directly applicable.  
\footnote{Alternatively, we may extend the $p_1$ model to directed networks with multiple communities.  However,   it is unclear 
what are appropriate models and approaches. In fact, due to the nonlinear $K_{ij}$, the $p_1$ model is hard to analyze  (e.g., \cite{JW2025}).  An extended version of the $p_1$ model will be even harder to analyze. \\}

The original weblog data \cite{weblog} is a directed
network  consisting of $1,490$ nodes (each node is a weblog) and $19,090$ directed edges (each edge is a hyperlink  representing).  Each node is manually labeled as democrat or republican, and therefore there are two perceivable communities: democrat and republican. 
Using the manual labels, we divide the network into two sub-networks, denoted by
Weblog-Dem and Weblog-Rep. As both sub-networks are directed networks 
with only $1$ community, it is reasonable to model them with the $p_1$ model. 
The citation network of statisticians is constructed using the recent MADStat data set
collected by Ji, Jin, Ke, and Li \cite{JBES}.  The MADStat contains the bibtex and citation
information of $83,331$ papers published in $36$ journals in $1971$-$2015$.
Ji {\it et al}. \cite{JBES} reported a citee network (which is undirected due to preprocessing)
with $2,831$ nodes (each is an author). They argued that (a) the $2,831$ nodes are the most active authors in their data range, (b) the network has $3$ interpretable communities: ``Bayes",  ``biostatistics" , and ``non-parametric".   
As the $p_1$ model is for directed networks with $1$ community, following \cite{JBES}, we first divide the dataset into these $3$ communities by applying the $k$-means
to the rows of the estimated membership matrix by Mixed-SCORE \cite{MSCORE}; see details therein.  
Despite that model in \cite{JBES} is undirected, we believe that the community structure there is similar to what we have here:  both networks contain three interpretable communities (Bayesian, biostatistics,  non-parametric) and an author in the Bayesian community (say) in one network is likely to be an author in the Bayesian community in the other.  Using these communities and original citation data in MADStat, we construct
six directed networks with two time period (i.e., 2013-2015, 2014-2015), denoted by Citation-Bay(13-15), Citation-Bio(13-15), and Citation-NP(13-15),
Citation-Bay(14-15), Citation-Bio(14-15), and Citation-NP(14-15), where a directed edge from $i$ to $j$ means that
papers published by author $i$ during 2013-2015 cite papers by author $j$ in the same period once (similar definition for 2014-2015 data). 
For reasons above, each network  is approximately a directed network with $1$ community so we can model it with the $p_1$ model (a further justification is desirable but very challenging; we leave it to the future). 
We apply our approaches to  all these networks, and the results are in Table \ref{tab:real}
(both $\hat{\sigma}^2 \equiv \widehat{V}_n(\rho_0)$ and the $p$-value are associated with the
hypothesis testing of $H_0: \rho = 0$, or $\rho_0 = 0$ by our previous notation).

\begin{table}
\centering
\caption{Results for real data analysis $(\hat{\sigma}^2 = \widehat{V}_n(\hat{\rho}))$.}
\label{tab:real}
\scalebox{1}{
\begin{tabular}{cccccccc}
Name  &  $\#$ of nodes & $\#$ of double  edges    & $\hat{\rho}^*$  &  $\hat{\sigma}$  &  $\psi_n^*$ &  $p$-value   \\
\hline
Weblog-Dem                 &   $ 758 $&$ 2177 $      &  $2.42 $           &$ 0.074 $&$ 32.5 $&$ < 10^{-4} $  \\
Weblog-Rep                 &   $ 732 $&$ 2244 $&  $2.74$           &$ 0.084 $&$ 32.5 $&$ < 10^{-4} $ \\
Citation-Bay(14-15)        & $ 919 $&$ 45 $&     $ 13.65 $&$ 0.607 $&$ 22.5 $&$ < 10^{-4}  $
\\
Citation-Bay(13-15)        &  $ 919 $&$ 96 $ &$ 4.64 $&$ 1.123 $&$ 4.1 $&$ < 10^{-4} $
\\
Citation-Bio(14-15)        & $ 1108 $&$ 33 $& $ 2.64 $&$ 1.604 $&$ 1.6 $&$ 0.100 $
\\
Citation-Bio(13-15)        & $ 1108 $&$ 98 $& $ 5.89 $&$ 1.118 $&$ 5.3 $&$ < 10^{-4}  $
\\
Citation-NP(14-15)         &  $ 804 $&$ 49 $ &$ 4.43 $&$ 1.414 $&$ 3.1 $&$ 0.002 $
   \\
Citation-NP(13-15)         &  $ 804 $&$ 118 $ & $ 3.37 $&$ 0.556 $&$ 6.1 $&$ < 10^{-4} $
\end{tabular}
}
\end{table}

For the two weblog sub-networks, the estimates of the reciprocity parameters, i.e., $\hat{\rho}^*$,   are $2.42$ and $2.74$, respectively. Also, for the null hypothesis $H_0: \rho=0$, both test statistics $\psi_n^*(\rho)$ for
two sub-networks are $32.5$, with very small $p$-values.  These results indicate a strong reciprocity.
Our results suggest the following. First, for  each of the networks Weblog-Dem and Weblog-Rep,
there is a strong reciprocity effect. The result is quite reasonable,
especially as the data set is collected in an election year, a time period where the bloggers with
similar political viewpoints try very hard to support each other. In such an atmosphere, it is very common that when one adds a hyperlink to another, he/she will get a hyperlink in return.
Despite that the weblog data has been well-studied in the literature (e.g., \cite{zhao2012,  JiangzhouWang2023Fast}),
our work is the first one that analyzes the reciprocal effects:
Many existing works treated the network as undirected by ignoring the directions of edges,
where the reciprocal effects are not properly modeled.  For citation networks, the $p$-values for
all $6$ networks are significant, except for Citation-Bio(14-15). This suggests that citations
 are indeed reciprocal in statistics, but it may take a long time (more than 1 year)  in order for the reciprocal effects to be fully realized. This is reasonable, for statistics is known to have a relatively slow publication cycle.
Also, note that the reciprocal effects in Bayes are generally stronger than those in
Non-parametric and Biostatistics, suggesting that Bayes community is more tightly-woven than the other two.

\section{Sketch of our proof ideas} 
\label{sec:proof}

To prove our results, the key is to bound $\mathbb{E}[(\hat{\rho} - \rho)^2]$ where $\hat{\rho}=\log(Q_n(a)/Q_n(b))$.  Note that $\hat{\rho} - \rho = \log(\frac{Q_n(a)}{e^{\rho} Q_n(b)}) = \log(1 + \frac{U_n(\rho)}{e^{\rho} Q_n(b)}) \approx 
 \frac{U_n(\rho)}{e^{\rho} Q_n(b)}$,  where $U_n(\rho) = Q_n(a) - e^{\rho} Q_n(b)$. 
Therefore, the key is to bound $\mathrm{Var}(Q_n(a))$, $\mathrm{Var}(Q_n(b))$ and $\mathrm{Cov}(Q_n(a), e^{\rho} Q_n(b))$. 
The following lemma is proved in the supplement  (recall that $(\mu, \nu, \eta)$ are defined in (\ref{plowrankmodel1})).

\begin{lemma} \label{lemma:snr2}
Suppose conditions (\ref{condition1}), (\ref{condition2}) and (\ref{condition4}) hold.  As $n \goto \infty$,
\begin{itemize}
\setlength \itemsep{-.5 em}
\item $\mathrm{Var}(Q_n(a))  = (1+o(1))\sum_{i \neq j} (2 r_{ij}^2 \Omega_{ij}^{11} + s_{ij}^2 \Omega_{ij}^{10})
=  (1+o(1)) [2 e^{\rho}\|\eta\|_1^2\|\mu\|_1^2\|\nu\|_1^2 +  e^{2\rho}\|\eta\|_1^4\|\mu\|_1\|\nu\|_1]$.
\item $\mathrm{Var}({Q}_n(b)) =    (1+o(1))  \sum_{i \neq j} t_{ij}^2 \Omega_{ij}^{10}
=  (1+o(1))  [(\mu,\eta)(\nu,\eta)\|\mu\|_1^2\|\nu\|_1^2 +3(\mu,\eta)\|\eta\|_1^2\|\mu\|_1\|\nu\|_1^2
 +3(\nu,\eta)\|\eta\|_1^2\|\mu\|_1^2\|\nu\|_1 +2\|\eta\|_1^4\|\mu\|_1\|\nu\|_1] + rem$, where the reminder term $|rem| \sim  \|\eta\|_3^3 \|\mu\|_1^2 \|\nu\|_1^2 = o(1) \cdot [e^{-2 \rho} \mathrm{Var}(Q_n(a)) + Q_n(b)]$.
 \item $\mathrm{Cov}(Q_n(a), e^\rho Q_n(b) )  =  - 2(1+o(1))\sum_{i\neq j} e^\rho s_{ij} t_{ij}\Omega_{ij}^{10} + o(1)\cdot [\mathrm{Var}(Q_n(a)) + e^{2\rho} \mathrm{Var}({Q}_n(b))]$.
\end{itemize}
\end{lemma}  
In the second item, the reminder term is $o(1) \cdot [e^{-2 \rho} \mathrm{var}(Q_n(a)) + \mathrm{Var}(Q_n(b))]$ so it has a negligible contribution in $\mathrm{Var}(U_n(\rho))$.  If we impose  
a mild condition of $a_n  \goto 0$, where $a_n: = \|\eta\|_3^3 / [(\mu, \eta)(\nu, \eta)] < 1$, 
then the reminder term is $o(1) \cdot \mathrm{Var}(Q_n(b))$.  
Lemma \ref{lemma:snr2} is the key to our proofs.  Take Theorem \ref{thm:snr1} for example: the 
key is to derive an exact formula for $\mathrm{Var}(U_n(\rho))$. But since $U_n(\rho) = Q_n(b) - e^{\rho} Q_n(a)$, 
the  formula follows directly from Lemma \ref{lemma:snr2}.  
    
\begin{table}[htb!]
	\renewcommand\arraystretch{0.9}
	\centering
	\resizebox{1.\textwidth}{!}{
		\begin{tabular}{|l|ll||l|ll|}
			\hline
			& Expression  & Variance  &  & Expression &   Variance\\
			\hline
			$S_{1,1}$ & $\sum_{i,j,k,\ell(dist)}W_{ij}^{11}W_{jk}^{00}W_{k\ell}^{01}W_{\ell i}^{00}$ &   $o(e^{\rho}\|\eta\|_1^2\|\mu\|_1^2\|\nu\|_1^2)$ &$R_{1,1}$ &  $\sum_{i,j,k,\ell(dist)}W_{ij}^{10}W_{jk}^{10}W_{k\ell}^{00}W_{\ell i}^{10}$ &   $o((\mu,\eta)(\nu,\eta)\|\mu\|_1^2\|\nu\|_1^2)$  \\
			$S_{1,2}$ & $\sum_{i,j,k,\ell(dist)}\Omega_{ij}^{11}W_{jk}^{00}W_{k\ell}^{01}W_{\ell i}^{00}$ & $o(e^{\rho}\|\eta\|_1^2\|\mu\|_1^2\|\nu\|_1^2)$ &$R_{1,2}$  &  $\sum_{i,j,k,\ell(dist)}\Omega_{ij}^{10}W_{jk}^{10}W_{k\ell}^{00}W_{\ell i}^{10}$  & $o((\mu,\eta)(\nu,\eta)\|\mu\|_1^2\|\nu\|_1^2)$\\
			$S_{1,3}$ & $\sum_{i,j,k,\ell(dist)}W_{ij}^{11}\Omega_{jk}^{00}W_{k\ell}^{01}W_{\ell i}^{00}$  & $o(e^{\rho}\|\eta\|_1^2\|\mu\|_1^2\|\nu\|_1^2)$ & $R_{1,3}$ & $\sum_{i,j,k,\ell(dist)}W_{ij}^{10}\Omega_{jk}^{10}W_{k\ell}^{00}W_{\ell i}^{10}$  & $o((\mu,\eta)(\nu,\eta)\|\mu\|_1^2\|\nu\|_1^2)$\\
			$S_{1,4}$ & $\sum_{i,j,k,\ell(dist)}W_{ij}^{11}W_{jk}^{00}\Omega_{k\ell}^{01}W_{\ell i}^{00}$ &   $o(e^{\rho}\|\eta\|_1^2\|\mu\|_1^2\|\nu\|_1^2)$ & $R_{1,4}$ &  $\sum_{i,j,k,\ell(dist)}W_{ij}^{10}W_{jk}^{10}\Omega_{k\ell}^{00}W_{\ell i}^{10}$  &   $o((\mu,\eta)(\nu,\eta)\|\mu\|_1^2\|\nu\|_1^2)$  \\
			$S_{1,5}$ & $\sum_{i,j,k,\ell(dist)}W_{ij}^{11}W_{jk}^{00}W_{k\ell}^{01}\Omega_{\ell i}^{00} $  &  $o(e^{\rho}\|\eta\|_1^2\|\mu\|_1^2\|\nu\|_1^2)$ & $R_{1,5}$ &  $\sum_{i,j,k,\ell(dist)}W_{ij}^{10}W_{jk}^{10}W_{k\ell}^{00}\Omega_{\ell i}^{10} $  &  $o((\mu,\eta)(\nu,\eta)\|\mu\|_1^2\|\nu\|_1^2)$  \\
			$S_{1,6}$ & $\sum_{i,j,k,\ell(dist)}\Omega_{ij}^{11}\Omega_{jk}^{00}W_{k\ell}^{01}W_{\ell i}^{00}$ &  $o(e^{2\rho}\|\eta\|_1^4\|\mu\|_1\|\nu\|_1)$ & $R_{1,6}$  & $\sum_{i,j,k,\ell(dist)}\Omega_{ij}^{10}\Omega_{jk}^{10}W_{k\ell}^{00}W_{\ell i}^{10}$  & $o((\mu,\eta)(\nu,\eta)\|\mu\|_1^2\|\nu\|_1^2)$\\
			$S_{1,7}$ & $\sum_{i,j,k,\ell(dist)}\Omega_{ij}^{11}W_{jk}^{00}\Omega_{k\ell}^{01}W_{\ell i}^{00}$ & $o(e^{2\rho}\|\eta\|_1^4\|\mu\|_1\|\nu\|_1)$ & $R_{1,7}$ &  $\sum_{i,j,k,\ell(dist)}\Omega_{ij}^{10}W_{jk}^{10}\Omega_{k\ell}^{00}W_{\ell i}^{10}$  & $o((\mu,\eta)(\nu,\eta)\|\mu\|_1^2\|\nu\|_1^2)$ \\
			$S_{1,8}$ & $\sum_{i,j,k,\ell(dist)}\Omega_{ij}^{11}W_{jk}^{00}W_{k\ell}^{01}\Omega_{\ell i}^{00}$  & $o(e^{2\rho}\|\eta\|_1^4\|\mu\|_1\|\nu\|_1)$ & $R_{1,8}$ & $\sum_{i,j,k,\ell(dist)}\Omega_{ij}^{10}W_{jk}^{10}W_{k\ell}^{00}\Omega_{\ell i}^{10}$  & $o((\mu,\eta)(\nu,\eta)\|\mu\|_1^2\|\nu\|_1^2)$\\
			$S_{1,9}$ & $\sum_{i,j,k,\ell(dist)}W_{ij}^{11}\Omega_{jk}^{00}\Omega_{k\ell}^{01}W_{\ell i}^{00}$  & $o(e^{\rho}\|\eta\|_1^2\|\mu\|_1^2\|\nu\|_1^2)$ &$R_{1,9}$  & $\sum_{i,j,k,\ell(dist)}W_{ij}^{10}\Omega_{jk}^{10}\Omega_{k\ell}^{00}W_{\ell i}^{10}$ & $o((\mu,\eta)(\nu,\eta)\|\mu\|_1^2\|\nu\|_1^2)$\\
			$S_{1,10}$ & $\sum_{i,j,k,\ell(dist)}W_{ij}^{11}\Omega_{jk}^{00}W_{k\ell}^{01}\Omega_{\ell i}^{00}$  & $o(e^{\rho}\|\eta\|_1^2\|\mu\|_1^2\|\nu\|_1^2)$ & $R_{1,10}$ & $\sum_{i,j,k,\ell(dist)}W_{ij}^{10}\Omega_{jk}^{10}W_{k\ell}^{00}\Omega_{\ell i}^{10}$  & $o((\mu,\eta)(\nu,\eta)\|\mu\|_1^2\|\nu\|_1^2)$\\
			$S_{1,11}$ & $\sum_{i,j,k,\ell(dist)}W_{ij}^{11}W_{jk}^{00}\Omega_{k\ell}^{01}\Omega_{\ell i}^{00}$ & $o(e^{\rho}\|\eta\|_1^2\|\mu\|_1^2\|\nu\|_1^2)$ & $R_{1,11}$  & $\sum_{i,j,k,\ell(dist)}W_{ij}^{10}W_{jk}^{10}\Omega_{k\ell}^{00}\Omega_{\ell i}^{10}$  & $o((\mu,\eta)(\nu,\eta)\|\mu\|_1^2\|\nu\|_1^2)$\\
			\multirow{2}{*}{$S_2$} & \multirow{2}{*}{$\sum_{i,j,k,\ell(dist)}\Omega_{ij}^{11}\Omega_{jk}^{00}\Omega_{k\ell}^{01}W_{\ell i}^{00}$} & $o(e^{\rho}\|\eta\|_1^2\|\mu\|_1^2\|\nu\|_1^2$
		  &\multirow{2}{*}{$R_2$}  & \multirow{2}{*}{$\sum_{i,j,k,\ell(dist)}\Omega_{ij}^{10}\Omega_{jk}^{10}\Omega_{k\ell}^{00}W_{\ell i}^{10}$}  & \multirow{2}{*}{$(1+o(1)) (\mu,\eta)\|\eta\|_1^2\|\mu\|_1\|\nu\|_1^2$}\\
		   & & $+e^{2\rho}\|\eta\|_1^4\|\mu\|_1\|\nu\|_1)$ & & &\\
			$S_3$ & $\sum_{i,j,k,\ell(dist)}\Omega_{ij}^{11}\Omega_{jk}^{00}W_{k\ell}^{01}\Omega_{\ell i}^{00}$  & $(1+o(1)) e^{2\rho}\|\eta\|_1^4\|\mu\|_1\|\nu\|_1$ & $R_3$  &  $\sum_{i,j,k,\ell(dist)}\Omega_{ij}^{10}\Omega_{jk}^{10}W_{k\ell}^{00}\Omega_{\ell i}^{10}$  & $o((\mu,\eta)(\nu,\eta)\|\mu\|_1^2\|\nu\|_1^2)$\\
			\multirow{2}{*}{$S_4$} &  \multirow{2}{*}{$\sum_{i,j,k,\ell(dist)}\Omega_{ij}^{11}W_{jk}^{00}\Omega_{k\ell}^{01}\Omega_{\ell i}^{00}$} & $o(e^{\rho}\|\eta\|_1^2\|\mu\|_1^2\|\nu\|_1^2$
			& \multirow{2}{*}{$R_4$}  & \multirow{2}{*}{$\sum_{i,j,k,\ell(dist)}\Omega_{ij}^{10}W_{jk}^{10}\Omega_{k\ell}^{00}\Omega_{\ell i}^{10}$}  & \multirow{2}{*}{$(1+o(1)) (\nu,\eta)\|\eta\|_1^2\|\mu\|_1^2\|\nu\|_1$}\\
			& & $+e^{2\rho}\|\eta\|_1^4\|\mu\|_1\|\nu\|_1)$ & & & \\
			\multirow{2}{*}{$S_5$} &   \multirow{2}{*}{$\sum_{i,j,k,\ell(dist)}W_{ij}^{11}\Omega_{jk}^{00}\Omega_{k\ell}^{01}\Omega_{\ell i}^{00}$} &  \multirow{2}{*}{$2(1+o(1)) e^{\rho}\|\eta\|_1^2\|\mu\|_1^2\|\nu\|_1^2$} & \multirow{2}{*}{$R_5$}  &  \multirow{2}{*}{$\sum_{i,j,k,\ell(dist)}W_{ij}^{10}\Omega_{jk}^{10}\Omega_{k\ell}^{00}\Omega_{\ell i}^{10}$}  & $   (1+o(1)) (\mu,\eta)(\nu,\eta) \|\mu\|_1^2 \|\nu\|_1^2 $\\
			& & & & & $ - (1 + o(1)) \|\eta\|_3^3  \|\mu\|_1^2 \|\nu\|_1^2$\\
			\hline
		\end{tabular}
	}
\caption{The variances of all $15$ terms in $Q_n(a)-\mathbb{E}[Q_n(a)]$ (left) and   $Q_n(b)-\mathbb{E}[Q_n(b)]$ (right).}\label{15-variances}
\label{tab:var}
\end{table}

We now sketch the high-level idea for proving Lemma \ref{lemma:snr2}, with
detailed proof in the supplement.   Consider $\mathrm{Var}(Q_n(a))$ (the first item of Lemma \ref{lemma:snr2}) first. Write
$Q_n(a) - \mathbb{E}[Q_n(a)] = S_1 + S_2 + S_3 + S_4 + S_{5}$,  where $S_1 = \sum_{k = 1}^{11} S_{1, k}$,  
and so $(Q_n(a) - \mathrm{E}[Q_n(a)])$ is  the sum of a  total of $15$ terms.
See Table \ref{tab:var} (left).   Here, $S_2, S_3, S_4, S_5$ are so-called {\it linear terms}.
Take $S_2 = \sum_{i, j, k, \ell (dist)} \Omega_{ij}^{11}  \Omega_{jk}^{00}  \Omega_{k\ell}^{01} W_{\ell i}^{00}$ for example.
Each terms in the big sum is a product of $4$ items: one is a $W$-term which is random, the other three are $\Omega$-terms which are non-random. Therefore, $S_2$ is a linear combination of many $W$-terms, and thus is linear.
In this sense, $S_{1,1}, S_{1,2}, \ldots, S_{1, 11}$ (and so $S_1$)  are {\it nonlinear terms}.
For each of these $15$ terms, we compute the variance separately; see Table \ref{tab:var}.
By delicate and lengthy analysis, we find that out of all $15$ terms, $S_3$ and $S_5$ are dominating in terms of variance,
where the correlation between $S_3$ and $S_4$ is asymptotically vanishing.
As a result,$\mathrm{Var}(Q_n(a)) = (1 + o(1)) \cdot \mathrm{Var}(S_3 + S_5) = (1 + o(1)) \cdot (\mathrm{Var}(S_3) + \mathrm{Var}(S_5))$.  
This largely reduces the difficulty level of downstream analysis. Finally, in our notations
$S_3 = \sum_{i \neq j} s_{ij} W_{ij}^{01}$ and $S_5 = \sum_{i \neq j} r_{ij} W_{ij}^{11}$. 
The first item of Lemma \ref{lemma:snr2} follows   by directly computing $\mathrm{Var}(\sum_{i \neq j} s_{ij} W_{ij}^{01})$, $\mathrm{Var}(\sum_{i \neq j} r_{ij} W_{ij}^{11})$, together
with basic calculus.  See Lemma D.1 in the supplement for the details.

Consider $\mathrm{Var}(Q_n(b))$ (the second item of Lemma \ref{lemma:snr2})  next.  Similarly, we write
$Q_n(b) - \mathbb{E}[Q_n(b)] = R_1 + R_2 + R_3 + R_4 + R_5$,  where $R_1 = \sum_{k=1}^{11} R_{1, k}$, 
 $R_2, R_3, R_4, R_5$ are linear terms and $R_{1, k}$ are nonlinear terms.
Similarly, we find that $(R_2 + R_4 + R_5)$  dominates in variance, and
$\mathrm{Var}(Q_n(b)) = (1 + o(1)) \cdot \mathrm{Var}(R_2 + R_4  + R_5)$.  
By our notations, it is seen
$R_2 + R_4 + R_5  = \sum_{i \neq j} t_{ij} W_{ij}^{10}$.  
The second result follows by directly computing $\mathrm{Var}(\sum_{i \neq j} t_{ij} W_{ij}^{10})$ and basic calculus.   See Lemma D.2  in the supplement for details.

We now consider $\mathrm{Cov}(Q_n(a)), e^{\rho} Q_n(b))$. By elementary statistics,
a term that has a negligible contribution in variance also has a negligible contribution in covariance. Therefore,
by the above discussions, $\mathrm{Cov}(Q_n(a),   e^{\rho} Q_n(b)) = \mathrm{Cov}(S_3 + S_5, e^{\rho} (R_2 + R_4 + R_5)) + o(1) \cdot (\mathrm{Var}(Q_n(a)) + e^{2 \rho} \mathrm{Var}(Q_n(b))$,  
where the first term on the right hand side is the leading term, and in our notations, it is
$\mathrm{Cov}(\sum_{i \neq j} (r_{ij} W_{ij}^{11} + s_{ij} W_{ij}^{10}), \sum_{i \neq j} t_{ij} W_{ij}^{10})$.  
The remaining part is detailed analysis. See Lemma D.3 in the supplement for details.

We now briefly explain why the proof of Lemma \ref{lemma:snr2} is so delicate and long.
First, we need to analyze $100+$ terms (each is a big sum similar as those in Table \ref{tab:var}) separately.
In detail, to compute $\mathrm{Var}(Q_n(a))$,  we have to split $(Q_n(a)- \E[Q_n(a)])$ into $15$ terms,  each is a big sum as in Table \ref{tab:var}, and analyze them separately; same for $\mathrm{Var}(Q_n(b))$.
Also, we also need analyze $\mathrm{Cov}(Q_n(a), e^{\rho} Q_n(b))$. Therefore,
for these together, we need to analyze at least $30+$ terms, each is a big sum similar to those in  Table \ref{tab:var}.
But this is not the end: for each of the $30+$ terms,
the variance further splits into several terms which we have to analyze  separately.
Take $\mathrm{Var}(S_{1,1}) =   \sum_{i,j,k,\ell(dist)}W_{ij}^{11}W_{jk}^{00}W_{k \ell }^{01}W_{\ell i}^{00}$ for example.  By the assumptions of the $p_1$ model,
$\mathrm{Var}(S_{1,1}) = \mathbb{E}[S_{1,1}^2]   = \mathbb{E}[ \sum_{i,j,k, \ell (dist)} \sum_{i', j', k', \ell' (dist)} W_{ij}^{11}W_{jk}^{00}W_{k \ell }^{01}W_{\ell i}^{00}W_{i'j'}^{11}W_{j'k'}^{00}\\
W_{k' \ell'}^{01}W_{\ell' i'}^{00}]$. 
To analyze the RHS, we need to use complicated combinatorics to further divide the RHS into several
sub-terms. In short, to prove Lemma \ref{lemma:snr2}, we need to carefully analyze $100+$ terms,  each  is a big sum similar to those in Table \ref{tab:var}.  Note that our statistics are U-statistics in nature, and
the above analytical approach is conventional: we are not aware of any shortcut.

Moreover,   for many of these terms,  we need to find the precise leading term (a cruder bound is frequently insufficient for the optimality results we desire). Also,
comparing the order of the variance of different terms is challenging. Take Table \ref{tab:var} (left) for example.
It is not immediately clear which of the $15$ listed variances are dominating.
Also, we face complicated combinatorics, especially as we do not have symmetry.  Take $S_{1,1} =  \sum_{i,j,k,\ell(dist)}W_{ij}^{11}W_{jk}^{00}W_{k \ell }^{01}W_{\ell i}^{00}$ for example, we do not have symmetry among the $4$ $W$ terms.  Last, from a practical viewpoint, we need to accommodate severe degree heterogeneity, so in our settings, different $(\alpha_i, \beta_i)$ may have very different magnitudes.
All these make our analysis hard, delicate, long, and error-prone.

However, this does not mean our proofs are hard to extend. For example,
for analysis of similar settings in the future, researchers can use a large part of our calculations.
Also, if we wish to sacrifice some of the assumptions, we may largely shorten the proofs. 

 {\bf Remark 8}. An interesting problem is how to derive an explicit approximation for $\mathrm{Var}(Q_n(a))$ and 
 	$\mathrm{Var}(Q_n(b))$ without assuming condition G1.  This is possible 
 	but the analysis is much more tedious. Take $\mathrm{Var}(Q_n(a))$ for example.  Without 
 	G1, we can only claim $\mathrm{Var}(Q_n(a)) = (1 + o(1)) \cdot [\sum_{k = 2}^5 \mathrm{Var}(S_k) + \sum_{2 \leq k \neq \ell \leq 5} \mathrm{Cov}(S_k, S_{\ell})]$. If we assume G1, then among the $10$ terms on the RHS, 
 	all except $\mathrm{Var}(S_3)$ and $\mathrm{Var}(S_5)$ are at the order of $o(1) \cdot \mathrm{Var}(Q_n(a))$. Therefore, $\mathrm{Var}(Q_n(a)) = (1 + o(1)) [\mathrm{Var}(S_3) + \mathrm{Var}(S_5)]$.



\section{Discussion} \label{sec:Discu}
Despite that it is one of the most popular models for directed network,  the $p_1$ model
does not model community structures. To overcome such a limitation,  we may
consider similar model as in (\ref{plowrankmodel1}), but $\widetilde{\Omega}$ is a general low-rank matrix which may be
much broader that that in (\ref{plowrankmodel1}).
Seemingly, the model is not only a generalization of the $p_1$ model but also the
popular block-models (e.g., \cite{SCORE-Review}).  In fact,
in  the special case where $K$ is the matrix of all ones, $\Omega$  is low-rank, and the model includes the block-model family 
as special cases. How to analyze such a model is of major interest. One
possible approach is to combine LCR with the idea of community detection,   where we first apply SCORE \cite{SCORE}  (say) for community detection, divide the nodes into different communities, and then apply the LCR for each community  (we may alternately apply community detection and LCR for several times). 
The $p_1$ model can also be viewed as a special case of the ERGM model \cite{Goldenberg2010,  ERGM}.
Despite its broadness,  the ERGM model is hard-to-analyze.  An interesting question is therefore to what extent our approach for the $p_1$ model is extendable to the ERGM model.  
Finally, our paper is also connected \cite{Feng2025reciprocitya, Feng2025reciprocityb}
(which appear when our paper is under review),  where \cite{Feng2025reciprocitya} studied a directed interaction model with sparsity and reciprocity, 
and \cite{Feng2025reciprocityb} studied reciprocity in the $p_1$ model with covariates.

\small 
\bibliographystyle{chicago}
\bibliography{reference}

%
%
%
%

\end{document}